\documentclass[pdflatex,sn-vancouver,Numbered]{sn-jnl}% Vancouver Reference Style

\usepackage{amsthm}%
\usepackage{mathrsfs}%
\usepackage[title]{appendix}%
\usepackage{textcomp}%
\usepackage{manyfoot}%
\usepackage{algorithmicx}%
\usepackage{algpseudocode}%
\usepackage{listings}%

\usepackage[utf8]{inputenc} % allow utf-8 input
\usepackage[T1]{fontenc}    % use 8-bit T1 fonts
\usepackage[table]{xcolor}

\usepackage{anyfontsize}
\usepackage{algorithm, algpseudocode, amsfonts, amsmath, amssymb, appendix}

\usepackage{bm, bbm, booktabs}
\usepackage{color, comment}
\usepackage{epstopdf, framed, graphicx, hyperref}
\definecolor{darkblue}{rgb}{0.10, 0.20, 0.65}
\definecolor{darkred}{rgb}{0.70, 0.00, 0.00}
\definecolor{darkgreen}{rgb}{0.20, 0.50, 0.20}
\hypersetup{
  colorlinks=true,
  breaklinks=true,
  linkcolor=darkred,
  urlcolor=darkblue,
  anchorcolor=darkblue,
  citecolor=darkblue
}
\usepackage{geometry}
\usepackage{latexsym, lineno}
\usepackage{microtype, multirow, multirow, nicefrac}
\usepackage{pdflscape, pifont}
\usepackage[figuresright]{rotating}
\usepackage{subfigure, pdflscape}
\usepackage{silence}
\usepackage{url, utfsym}

\DeclareMathOperator*{\argmin}{arg\rm{}min}

\newcommand{\bs}{\mathbf{s}}

\newcommand{\bu}{\boldsymbol{u}}

\newcommand{\bX}{\mathbf{X}}

\begin{document}

\title[LES-SINDy]{LES-SINDy: Laplace-Enhanced Sparse Identification of Nonlinear Dynamical Systems}

\author[1]{\fnm{Haoyang} \sur{Zheng}}\email{zheng528@purdue.edu}

\author*[1, 2]{\fnm{Guang} \sur{Lin}}\email{guanglin@purdue.edu}

\affil[1]{\orgdiv{School of Mechanical Engineering}, \orgname{Purdue University}, \orgaddress{\street{585 Purdue Mall}, \city{West Lafayette}, \postcode{47907}, \state{IN}, \country{USA}}}

\affil[2]{\orgdiv{Department of Mathematics}, \orgname{Purdue University}, \orgaddress{\street{150 North University Street}, \city{West Lafayette}, \postcode{47907}, \state{IN}, \country{USA}}}

%%==================================%%
%% Sample for unstructured abstract %%
%%==================================%%

\abstract{

Sparse Identification of Nonlinear Dynamical Systems (SINDy) is a powerful tool for data-driven discovery of governing equations. However, it encounters challenges when modeling complex dynamical systems involving high-order derivatives or discontinuities, particularly in the presence of noise. These limitations restrict its applicability across various fields in applied mathematics and physics. To mitigate these, we propose Laplace-Enhanced SparSe Identification of Nonlinear Dynamical Systems (LES-SINDy). By transforming time-series measurements from the time domain to the Laplace domain using the Laplace transform and integration by parts, LES-SINDy enables more accurate approximations of derivatives and discontinuous terms. It also effectively handles unbounded growth functions and accumulated numerical errors in the Laplace domain, thereby overcoming challenges in the identification process. The model evaluation process selects the most accurate and parsimonious dynamical systems from multiple candidates. Experimental results across diverse ordinary and partial differential equations show that LES-SINDy achieves superior robustness, accuracy, and parsimony compared to existing methods.
}


\maketitle

\section{Introduction}\label{sec_intro}

% \noindent\textbf{Model Discovery:}\ \ 
The discovery of scientific laws from measurements is a significant intellectual milestone, and its motivation arises from the widespread occurrence of nonlinear dynamical systems in science and engineering. Understanding the governing equations, which often take the form of ordinary differential equations (ODEs), partial differential equations (PDEs), and stochastic differential equations (SDEs), is essential for accurate prediction, effective control, and informed decision-making \citep{brunton2022data, kutz2013data}. In many complex systems, the underlying dynamics remain poorly understood, which renders conventional modeling techniques based on first principles both challenging and, at times, intractable. To tackle the challenge of model discovery in dynamical systems, Sparse Identification of Nonlinear Dynamics (SINDy) \citep{brunton2016discovering} offers a data-driven solution. By the use of given measurements, SINDy constructs parsimonious models that capture the essential features of system dynamics without the need for detailed knowledge of the underlying physics. The strength of SINDy lies in its ability to identify sparse and interpretable models, based on the assumption that the system's dynamics can be represented as a sparse linear combination of candidate functions. This process involves iterative optimization through sparse regression \citep{daubechies2004iterative} and the selection of the most relevant terms from a comprehensive library, which enables the discovery of governing equations that are both accurate and physically meaningful.

Building on the idea of using sparse regression techniques to discover nonlinear dynamical systems, extensive research has been conducted to enhance the SINDy framework for various objectives or to apply it across diverse domains. Schaeffer \cite{schaeffer2017learning} employed lasso regression within SINDy to identify PDEs. Rudy \textit{et al.} \cite{rudy2017data, rudy2019data} advanced the field by proposing PDE functional identification of nonlinear dynamics (PDE-FIND), which extends SINDy to discover dynamical systems from spatial time series measurements.  Brunton \textit{et al.} \cite{brunton2016sparse} and Kaiser \textit{el al.} \cite{kaiser2018sparse} focused on generalizing SINDy to include control inputs and facilitating its use in model predictive control. Zhang \textit{et al.} \cite{zhang2018robust} and Kaheman \textit{et al.} \cite{kaheman2020sindy} extended SINDy to identify implicit differential equations (DEs), thereby broadening its applicability. When modeling implicit DEs, criteria such as the Akaike information criterion (AIC) \cite{akaike1974new, schwarz1978estimating} and Bayesian information criterion (BIC) \cite{akaike1998information, mangan2017model} are commonly employed to select the most accurate and parsimonious model among multiple candidates. To extract dynamical systems from complex, high-dimensional measurements, methods like autoencoders \cite{champion2019data} and physics-informed neural networks \cite{raissi2019physics, raissi2018deep, chen2021physics} have effectively employed automatic differentiation to approximate system variables and identify significant derivative terms. More recent research has integrated neural networks with manifold theory to directly discover intrinsic state variables and their dynamics from time series measurements \cite{floryan2022data}. To enhance the utility of the SINDy framework in various scientific and engineering disciplines, the Python package PySINDy \cite{desilva2020, Kaptanoglu2022} was developed to provide flexible tools for data-driven model discovery. With these advanced tools, data-driven methods can be applied effectively in diverse domains, including control systems \cite{markovsky2021behavioral, fasel2021sindy}, fluid mechanics \cite{brunton2020machine}, biomechanics \cite{zheng2022data, liang2022discovering}, and chemistry \cite{keith2021combining}.

\vspace{0.05 in}
\noindent\textbf{Challenges and Motivations:}\ \ While the SINDy framework is adaptable to various applications, ensuring robustness and accuracy is still a concern, particularly in the presence of high-order derivatives and discontinuities in the system dynamics \cite{rudy2017data}. The implementation of SINDy assumes that all terms in the library are smooth and differentiable. However, this assumption presents challenges when numerical differentiation methods, such as finite differences, are used to approximate high-order derivatives or discontinuous terms in the library. Consequently, identifying governing equations that include such terms becomes difficult, as such approximate errors can significantly impact the accuracy of the resulting sparse model. Furthermore, the presence of noisy measurements further exacerbates these issues, which potentially leads to instability and inaccuracy in the approximations. 

To address these challenges, researchers have explored various strategies, and one is to learn multiple SINDy models and effectively leverage them. For example, Ensemble-SINDy \cite{fasel2022ensemble} discovered multiple SINDy models from different subsets of measurements simultaneously, and then aggregated the results to form the final model. Additionally, Uncertainty Quantification SINDy (UQ-SINDy) \cite{hirsh2022sparsifying} introduced a probabilistic assessment of each candidate function, which was later incorporated with physical constraints \cite{zheng2024constrained}. By incorporating Bayesian methods to analyze multiple potential models, UQ-SINDy ensures more reliable model discovery and provides confidence intervals to quantify the certainty of the results. While these methods enhance robustness and reliability in model identification, they often become computationally intensive as the number of measurements accumulates. The computational cost is largely driven by the sparse regression process in SINDy, where the expense of matrix operations, like inversion, scales roughly cubicly with more measurements. Furthermore, the need for repeated sparse regression across numerous iterations or samples further amplifies the computational cost.  

An alternative strategy is to first transform the library from the time domain to another domain, such as the frequency or Laplace domain, where high-order derivatives and discontinuities are more easily handled. Sparse regression is then employed in the transformed domain. For instance, the Weak-SINDy \cite{messenger2021weak, messenger2022weak} transforms the SINDy library to the frequency domain using Fourier basis functions \cite{li2020fourier, li2022learning}, followed by sparse regression to produce the model. Intuitively. this approach is more efficient because the computational burden is no longer dictated by the sparse regression process. In the frequency domain, the matrix structure may be simplified and independent of the number of measurements (resolution-free), which makes matrix operation computationally advantageous.

Inspired by the successful use of Laplace transformations in different machine-learning tasks \cite{holt2022neural, cao2023lno}, we propose Laplace-Enhanced Sparse Identification of Nonlinear Dynamical Systems (LES-SINDy). This framework is designed to achieve robustness, efficiency, and accuracy in data-driven model discovery. Given time-series measurements, LES-SINDy first constructs candidate functions in the time domain and then subsequently transforms them to the Laplace domain using the Laplace transformation and integration by parts. After the transformation, a resolution-free sparse regression is employed to identify several potential candidate models. Lastly, the most appropriate model is selected from candidate models using a modified AIC.

It is important to note that both LES-SINDy and Weak-SINDy aim to identify dynamical systems by transforming the library from the time domain to the frequency domain. However, we provide further insights into why LES-SINDy is a more effective option for discovering nonlinear dynamical systems, which will be further supported by the experimental results. One advantage of the Laplace transform over the Fourier transform is its ability to effectively handle functions that exhibit unbounded growth. The Fourier transform requires that all the candidate functions be absolutely integrable over the time domain, which might be problematic to approximate functions that grow exponentially or do not decay sufficiently. Consequently, there is no clear evidence to support Weak-SINDy can accurately identify dynamical systems with unbounded growth, such as \eqref{eq:homo_fourth} (Fig. \ref{fig:fourth}) and \eqref{eq:sinh_func}-\eqref{eq:cosh_func} (Fig. \ref{fig:hyperbolic}). In contrast, the Laplace transform, with its exponentially decaying weighting factor, can handle such candidate functions by selecting $s$ with an appropriately large real part. This makes the Laplace transform more suitable for analyzing candidate functions that do not satisfy the integrability conditions required by the Fourier transform.
% Another inherent advantage of the Laplace transformation lies in its weighting scheme, which is beneficial 
Moreover, the exponential weighting scheme in LES-SINDy can also benefit numerically simulated measurements, which are prone to accumulate errors over time. By assigning greater weight to function values at earlier times, the Laplace transform focuses on the more accurate, less error-prone parts of the simulation. This is in contrast to the Fourier transform, which may inadvertently introduce accumulated errors, particularly in the later stages of the simulation. Therefore, in scenarios where maintaining numerical stability is crucial, the Laplace transformation can be a more reliable choice.
% \subsection{Contributions}

\vspace{0.05 in}
\noindent\textbf{Contributions:}\ \ We summarize the major contributions to this work as follows.
\begin{enumerate}
    \item We propose the novel LES-SINDy framework for robust and accurate data-driven model discovery, which leverages the Laplace transform and integration by parts to transform the candidate functions from the time domain to the Laplace domain.
    \item LES-SINDy addresses key challenges in existing methods by handling high-order derivatives,  discontinuities, unbounded growth functions, and accumulated numerical errors.
    \item LES-SINDy is ``resolution-free'': it is unaffected by measurement resolution, while the computational burden in SINDy increases cubically with increasing measurements.
    \item We test LES-SINDy across a diverse range of dynamical systems to validate its efficiency and robustness, including high-order ODEs, ODEs with discontinuous inputs, ODEs with trigonometric and hyperbolic functions, nonlinear ODE systems, and PDEs.
    % \item \textcolor{red}{A Python library for LES-SINDy, compatible with PySINDy \cite{desilva2020, Kaptanoglu2022}, will be released on GitHub to support the research community after publication.}
\end{enumerate}

\section{Results}\label{sec_exp}

In this section, we conduct a series of experiments to evaluate the performance of LES-SINDy across a range of dynamical systems. These include higher-order ODEs, ODEs with discontinuous inputs, ODEs involving trigonometric and hyperbolic functions, nonlinear ODE systems, and PDEs. The experiments were run on a server featuring an Intel(R) Core(TM) i9-14900K processor and 128 GB DDR4 memory. 
% The detailed experimental setup is summarized in Table \ref{tab:exp_summary} in Appendix \ref{appendix:exp}. 
For more experimental details, interested readers can refer to Appendices from \ref{appendix:high-order} to \ref{appendix:pde}.

\subsection{General framework of LES-SINDy}

\begin{figure}[!htbp]
    \hspace*{-0.0cm}
    \includegraphics[width=1.0\linewidth]{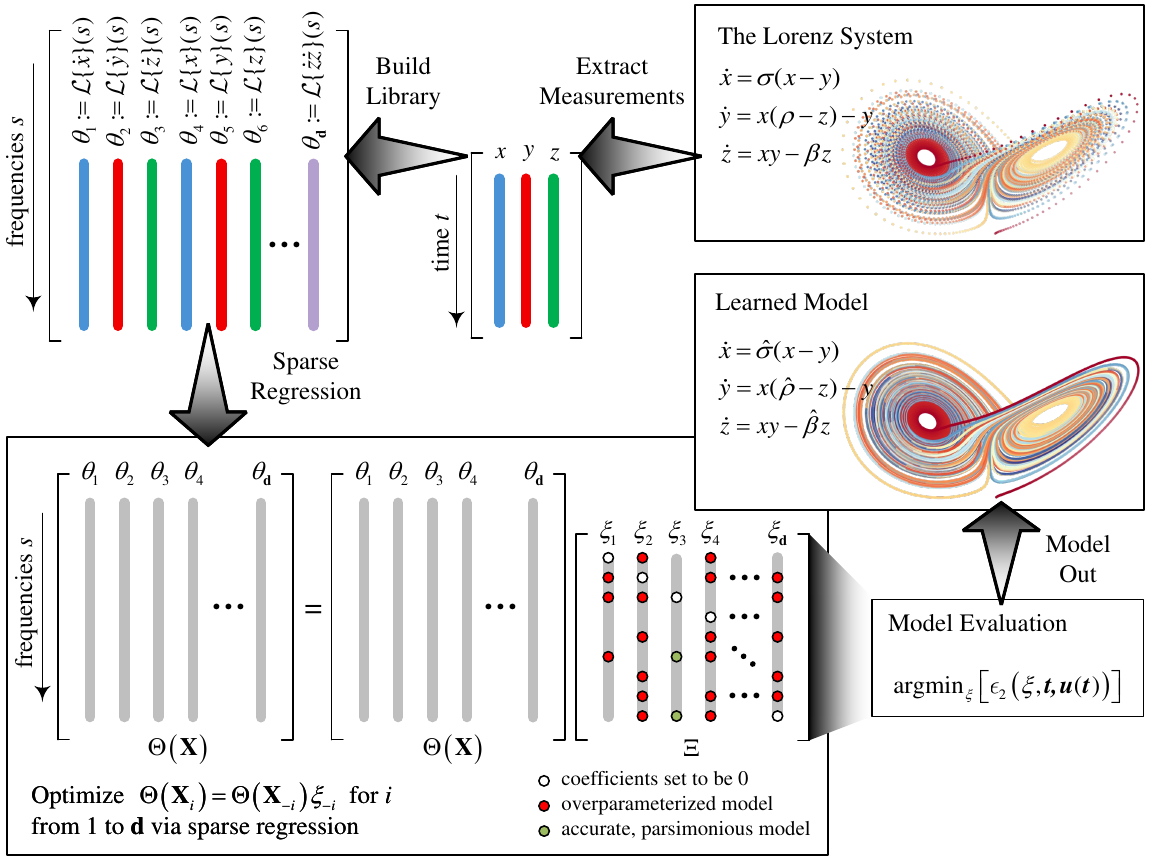}
    \caption{Illustration of the LES-SINDy framework, demonstrated on identifying the Lorenz system. The process starts by extracting time-series measurements from sensors, which are used to form a Laplace-enhanced library. Sparse regression techniques are then employed to identify potential models in the Laplace domain. The evaluation stage determines the most accurate and parsimonious model, which successfully recovers the Lorenz system.}
    \label{fig:framework}
\end{figure}
The proposed LES-SINDy framework aims to identify implicit DEs from time-series measurements. As shown in Fig. \ref{fig:framework}, LES-SINDy begins by constructing a comprehensive library of candidate functions derived from the measurements in the time domain, which encompasses state variables, their derivatives, and other relevant features that characterize the system dynamics. To further enhance the expressiveness of this library, these candidate functions are transformed into the Laplace domain, which provides an alternative representation that is particularly effective for capturing the dynamics of complex systems. Sparse regression is employed in the model discovery process to ensure that the final model captures the essential dynamics without unnecessary complexity. The iterative nature of sparse regression allows for systematic refinement of the model, which progressively eliminates less significant terms until an optimal balance between complexity and accuracy is achieved. The resulting models are evaluated by AIC with correction (AICc) based on both their accuracy and parsimony. 

\subsection{Experimental setup}

To evaluate the LES-SINDy framework, a series of experiments were conducted using a diverse range of DEs. The scope and details of these experiments are summarized in Table \ref{tab:exp_summary}, which demonstrates the framework's robustness and versatility across five different types of ODEs and PDEs: High-order ODEs were included to test LES-SINDy's capability to handle systems with multiple derivatives, which are often found in mechanical and physical systems. Discontinuous ODEs were motivated by the need to capture abrupt changes or impulses in the system, which are typical in control systems and signal processing. Nonlinear ODE systems were selected to assess LES-SINDy's performance on complex, interacting dynamical systems that exhibit chaotic or predator-prey dynamics. PDEs were used to evaluate the framework's ability to handle spatiotemporal dynamics, where both temporal and spatial derivatives play significant roles. These diverse experimental setups provide a comprehensive assessment of LES-SINDy's effectiveness across a wide range of dynamical systems and conditions.

\begin{table*}[!htbp]
\centering
% \hspace{-2mm}
\caption{Summary of Models and Experimental Setups for LES-SINDy Evaluation. This table outlines the ODE/PDE models (categorized by high-order ODEs, discontinuous ODEs, trigonometric and hyperbolic ODEs, nonlinear ODE systems, and PDEs) and experimental conditions (measurement noise, discretization schemes, etc.) used to assess the performance of the LES-SINDy algorithm, which covers a wide range of dynamical systems and conditions.}\label{tab:exp_summary}
\resizebox{1.05\textwidth}{!}{
\begin{tabular}{cllccc}
\toprule
type    &     model name    & form & noise & temporal grid & spatial grid \\ \midrule
\multirow{2}{*}{high-order ODEs}             
    & Duffing oscillator        & ${u}_{tt} + \alpha {u}_{t} + \beta u = 0$ &  $20\%$ &  $[0, 10]$\ \ $m=1000$ & $\usym{1F5F6}$ \\
    & fourth-order ODEs             & ${u}_{tttt} + \alpha {u}_{tt} + \beta u = 0$ &  clean  &  $[0, 20]$\ \ $m=200$ & $\usym{1F5F6}$ \\ \midrule
\multirow{2}{*}{\begin{tabular}[c]{@{}c@{}}ODEs with \\
 discontinuous functions\end{tabular}}           
    & delta              & ${u}_{tt} + \alpha {u}_{t} + \beta u - F_0 \delta(t_0) = 0$ & $10\%$ & $[0, 10]$\ \ $m=1000$ &  $\usym{1F5F6}$ \\
    & step               & ${u}_{t} + \alpha u - F_0 H(t_0) = 0$  & $20\%$ & $[0, 10]$\ \ $m=1000$ & $\usym{1F5F6}$ \\ \midrule
\multirow{4}{*}{\begin{tabular}[c]{@{}c@{}}ODEs with trigonometric\\ 
or hyperbolic functions\end{tabular}} 
    & sine               & ${u}_{tt} + 15u = 2\sin(3t)$ & $10\%$ &  $[0, 10]$\ \ $m=1000$ &  $\usym{1F5F6}$ \\
    & cosine             & ${u}_{tt} + 4u = \cos(t)$ & $10\%$ &  $[0, 10]$\ \ $m=1000$ & $\usym{1F5F6}$ \\
    & sinh               & ${u}_{tt} + 4u = \sinh(2t)$ & clean &  $[0, 100]$\ \ $m=10000$  & $\usym{1F5F6}$ \\
    & cosh               & ${u}_{tt} - 4u = \cosh(2t)$ & clean  & $[0, 100]$\ \ $m=10000$ & $\usym{1F5F6}$ \\ \midrule
\multirow{3}{*}{nonlinear ODE systems}  
    & Lorenz             & \begin{tabular}[c]{@{}ll@{}}$x_t = \sigma (y - x),$\\ $y_t = x (\rho - z) - y,$\\ $z_t = x y - \beta z.$\end{tabular} & {clean} & $[0, 100]$\ \ $m=10000$ & $\usym{1F5F6}$ \\ \cmidrule{2-6}
    & Lotka-Volterra      & \begin{tabular}[c]{@{}ll@{}}$x_t = \alpha x - \beta xy,$\\ $y_t = \delta xy - \gamma y.$\end{tabular}     & {clean} & $[0, 100]$\ \ $m=10000$  &  $\usym{1F5F6}$ \\ \midrule
\multirow{4}{*}{PDEs}                                                          
    & convection-diffusion  & $\frac{\partial u}{\partial t} + c \frac{\partial u}{\partial x} - D \frac{\partial^2 u}{\partial x^2} = 0$ &  $30\%$ & $[0, 5]$\ \ $m=500$  &  $[0, 5]$\ \ $n=1000$  \\ [1.2ex]
    & Burgers              & $\frac{\partial u}{\partial t} + u \frac{\partial u}{\partial x} - \nu \frac{\partial^2 u}{\partial x^2} = 0$ & $20\%$ & $[0, 5]$\ \ $m=500$ &  $[0, 20]$\ \ $n=1000$  \\ [1.2ex]
    & Kuramoto-Sivashinsky & $\frac{\partial u}{\partial t} + u \frac{\partial u}{\partial x} + \frac{\partial^2 u}{\partial x^2} + \frac{\partial^4 u}{\partial x^4} = 0$ &  $20\%$ & $[0, 100]$\ \ $m=1024$ & $[0, 100]$\ \ $n=251$ \\  [1.0ex] \bottomrule
\end{tabular}
}
\end{table*}

\subsection{High-order ODEs}

Higher-order ODEs are frequently encountered in the modeling of physical systems, where the dynamics are governed by second-order, third-order, or even higher-order derivatives. These equations often describe complex phenomena such as mechanical vibrations, wave propagation, and electrical circuits. We focus on two representative examples to evaluate the performance of LES-SINDy in identifying higher-order ODEs.

\begin{figure*}[!ht]
\centering
\subfigure[The Duffing oscillators with different $\zeta$ values.]{
\centering
\includegraphics[height=0.22\columnwidth]{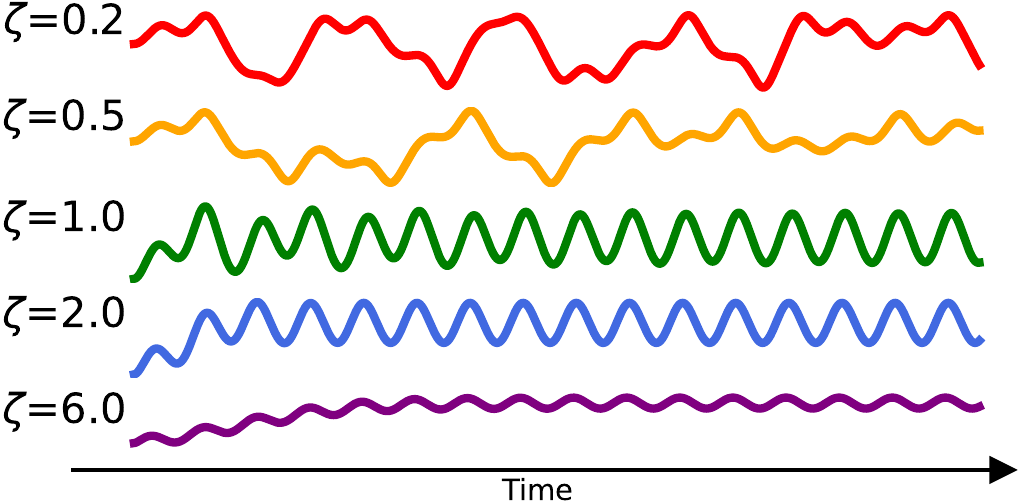}\label{fig:duffing}
}\ \ \ \ \ \ 
\subfigure[The fourth-order linear ODE.]
{
\centering
\includegraphics[height=0.22\columnwidth]{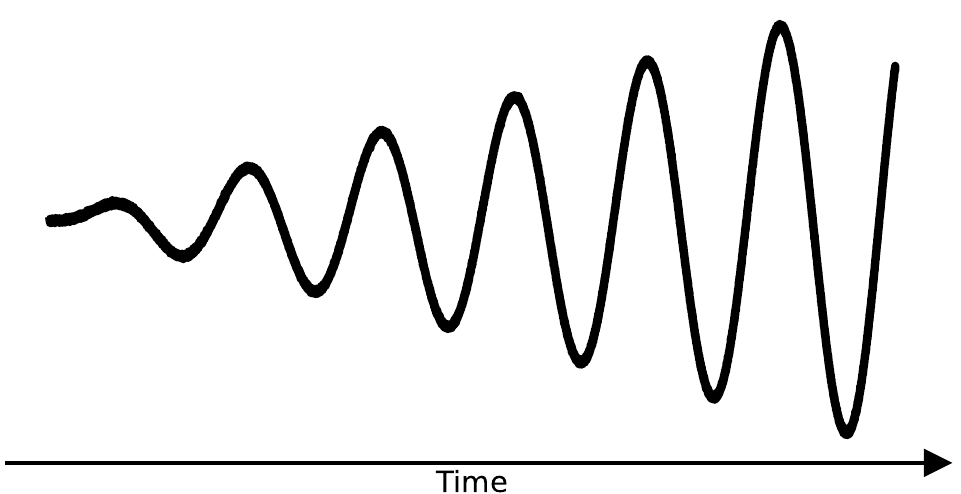}\label{fig:fourth}
}

\subfigure[log RMSE error mappings.]
{
\centering
\includegraphics[height=0.35\columnwidth]{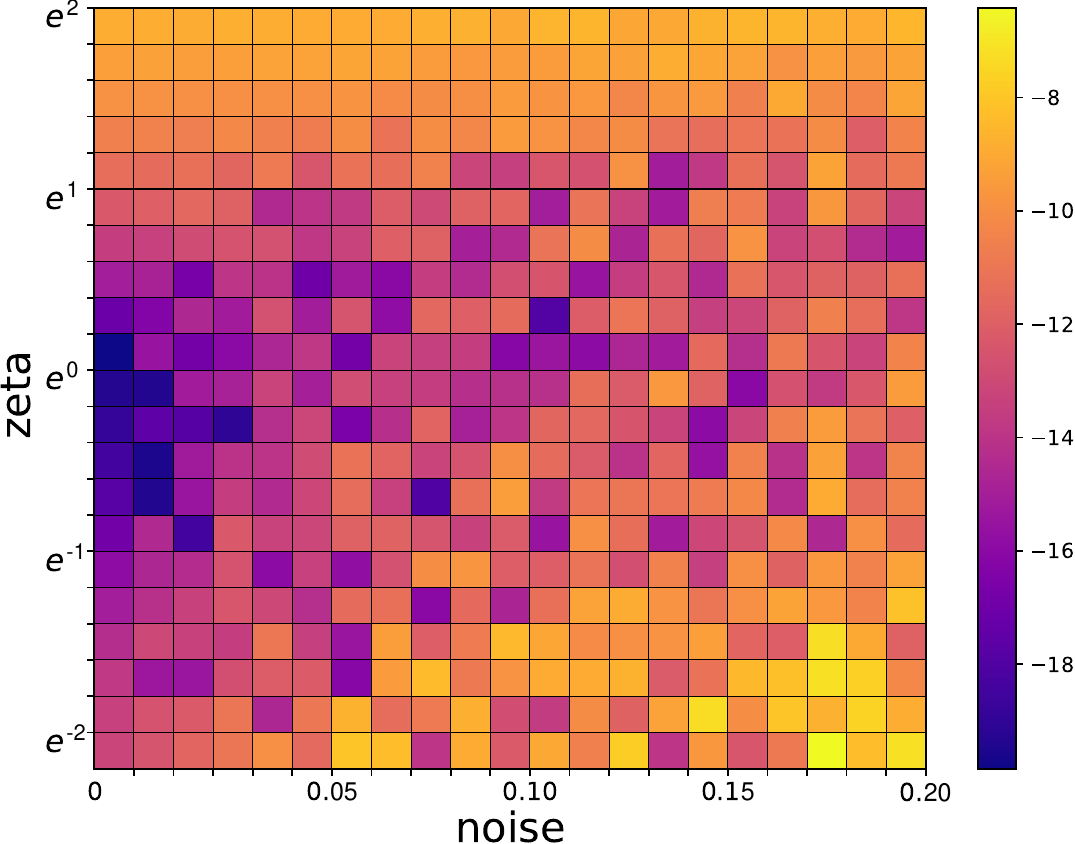}\label{fig:duffing_zeta}\ \ \ \ \ \ 
\includegraphics[height=0.35\columnwidth]{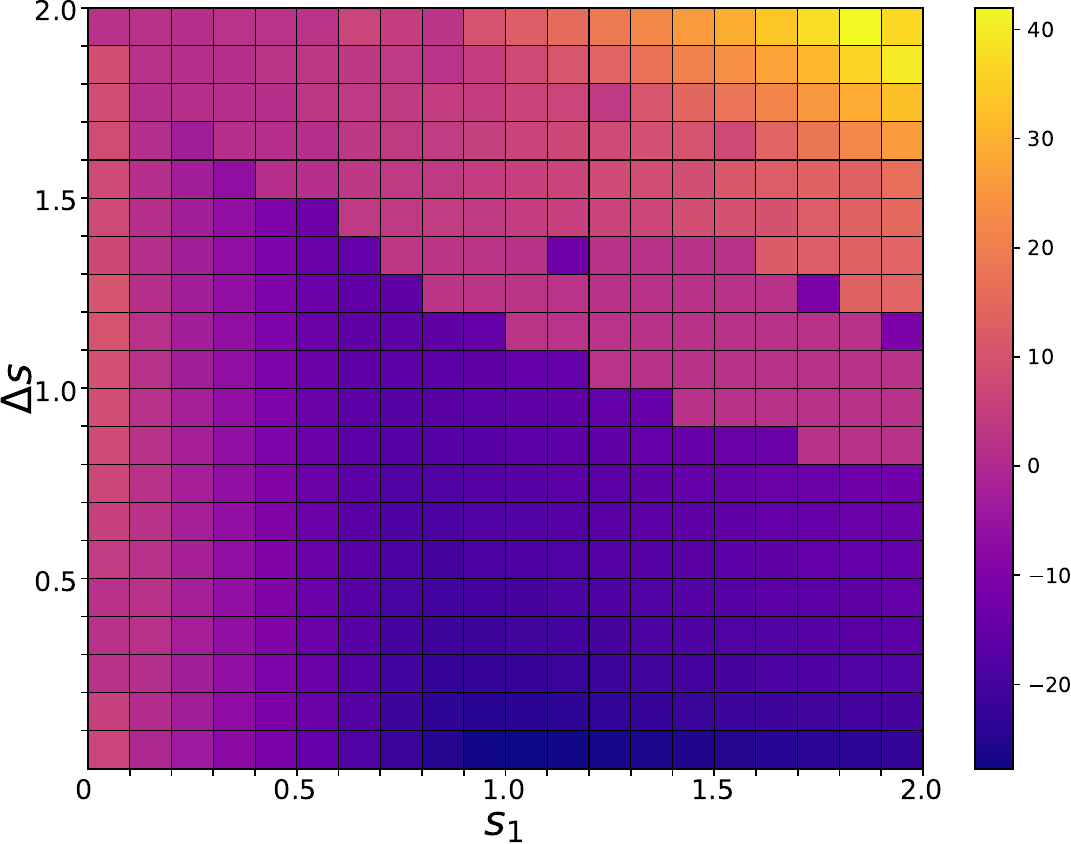}\label{fig:errorbar_fourth}
}
% \vspace{-0.15 in}
\caption{Simulation of the High-Order Linear ODE Systems and the corresponding log RMSE mappings. \textbf{Top Left:} Dynamic behavior of the Duffing oscillator across various damping ratios $\zeta$. The x-axis denotes the time, and the y-axis represents the magnitude of the dynamical system responses. When $\zeta=0.20$, it represents a lightly damped system, where oscillations will persist with gradual decay. $\zeta=0.50$ represents a moderately damped system, where the oscillations decay faster than in the lightly damped case. For $\zeta=1.0$, the system returns to equilibrium as quickly as possible without oscillating. Overdamped responses are demonstrated at $\zeta = 2.0$ and $\zeta = 6.0$, with increasingly slow returns to equilibrium, which indicates higher damping levels. \textbf{Top Right:} Dynamic behavior of a fourth-order linear ODE. The response implies quasi-periodic motion, which is characterized by a combination of oscillations and a time-dependent modulation of amplitude. This results in oscillations whose amplitude gradually increases over time, which reflects a deterministic and predictable pattern. The system exhibits intricate behavior where the oscillations shift in both amplitude and phase and demonstrates the dynamic interplay between periodicity and linear modulation as time progresses. \textbf{Bottom Left:} log RMSE Error Mapping for the Duffing Oscillators: Noise Scale and  $\zeta$ Value Effects in Logarithmic Scale. The x-axis denotes the selection of the noise scale, the y-axis is the magnitude of $\zeta$, and the colors represent the log RMSE errors between the predicted results and the ground truth. Magnitudes of noise scale increase uniformly from 0.0 to 0.20, and the $\zeta$ increases exponentially from $e^{-2}\approx 0.135$ to $e^2\approx 7.389$. The heat mapping indicates the $\zeta$ value may affect the result of good predictions. \textbf{Bottom Right:} {log RMSE Error Mapping for a Homogeneous Linear Fourth-Order ODEs: Initial Value and Step Size Effects in Logarithmic Scale. The x-axis denotes the selection of the initial value $s_1$ (intercept), the y-axis is the interval between different $s$ values (step size), and the colors represent the log RMSE errors between the predicted results and the ground truth. Values of $s$ are uniformly distributed from $s_1$ to $s_{20}$ based on the chosen intercept and the step size. The heat mapping indicates the appropriate selection of intercept and step size is crucial for good predictions.}}
\label{fig:high_order}\vspace{-0.0 in}
\end{figure*}

\textbf{The Duffing Oscillator} in Fig. \ref{fig:duffing}, as a well-known model in nonlinear dynamics, is widely used to describe systems with a nonlinear restoring force, such as certain types of mechanical oscillators. The second-order ODE governing the Duffing oscillator is typically expressed as:
\begin{equation}\label{eq:duffing}
    % {u}_{tt} + \gamma {u}_{t} + \alpha u + \beta u^3 = \delta \cos(\omega t),
    {u}_{tt} + 2\omega \zeta {u}_{t} + \omega^2 u + \alpha u^3 = \delta \cos(\omega_d t),
\end{equation}
where $u:=\mathbb R_{\geq 0}\rightarrow \mathbb R$ denotes the displacement at time $t$, $\omega$ is the system’s natural frequency, $2\omega \zeta$ represents the damping effect, $\alpha$ characterizes the stiffness, $\delta$ is the amplitude of the driving force, and $\omega_d$ is the driving frequency. 

In certain scenarios, especially when analyzing the system's behavior in the absence of an external driving force (i.e., when $\delta = 0$), the Duffing oscillator can be simplified to a form that is easier to handle both analytically and numerically. By assuming the absence of the nonlinear cubic term (i.e., $\alpha = 0$) and external forces, the equation can be rewritten as:
\begin{equation}\label{eq:duffing1}
    {u}_{tt} + 2\omega \zeta {u}_{t} + \omega^2 u = 0,
\end{equation}
where the system simplifies to ignore a linear restoring force and a damping force proportional to the velocity. 

It should be noted that the selection of the damping coefficient $\zeta$ plays a crucial role in determining the behavior of the system (see dynamical behaviors with different $\zeta$ values shown in Fig. \ref{fig:duffing}). When $\zeta$ is small, the system shows underdamped behavior, with oscillations that slowly diminish in amplitude. The frequency is slightly below the natural frequency $\omega$, and the energy loss is gradual. At $\zeta = 1.0$, the system is critically damped. It returns to its equilibrium position swiftly and without oscillations, which is desired in scenarios requiring quick stabilization. In contrast, large $\zeta$ values result in overdamping, where the system slowly reaches equilibrium without oscillations. This is typical in systems engineered to suppress vibrations.

One interesting observation is that the choice of $\zeta$ may also impact the performance of LES-SINDy. Following \eqref{eq:duffing1} and with $\omega = 2.0$, we varied the magnitude of $\zeta$ and observed the resulting prediction. The prediction is further compared with the measurements and evaluated with the logarithm of the root mean square error (log RMSEs). As shown in Fig. \ref{fig:duffing_zeta} (left), LES-SINDy performs best when $\zeta$ is close to $e^0=1.0$. However, too small or large $\zeta$ values lead to deviations from the true measurements. This occurs because LES-SINDy uses thresholding least squares for sparsification, which may incorrectly eliminate important coefficients when $\zeta$ is not comparable to other parameters.

\textbf{Homogeneous Linear Fourth-Order ODEs} (Fig. \ref{fig:fourth}) are essential in various engineering and physics applications, particularly in the study of beam deflections, vibrations, and stability analysis. We consider the following fourth-order ODE as a test case:
\begin{equation}\label{eq:homo_fourth}
    {u}_{tttt} + \alpha {u}_{tt} + \beta u = 0.
\end{equation} 

This homogeneous linear equation represents a system with complex characteristic roots, which leads to solutions that exhibit oscillatory behavior with specific decay or growth rates. Specifically, the coefficient $\alpha$ corresponds to the second derivative term ${u}_{tt}$, which is commonly related to velocity or inertia in physical systems. $\alpha$ affects the system's damping characteristics: a larger $\alpha$ increases damping, which leads to quicker stabilization with diminished oscillations, while a smaller $\alpha$ permits greater oscillatory motion due to lower energy dissipation. The parameter $\beta$ is associated with the $u$ term, which often represents the system’s restoring force or stiffness. In mechanical systems, $\beta$ correlates with the rigidity or stiffness of the structure. A higher $\beta$ typically results in a stiffer, higher-frequency response, whereas a lower $\beta$ produces slower oscillations with a more relaxed system response. 

\begin{table}[!htbp]
\centering\caption{Results of LES-SINDy applied to a high-order linear ODE.}\label{tab:high_order}
    \begin{tabular}{l|cccccc}
        \toprule
        \makebox[0.10\textwidth][l]{results}     &   \makebox[0.05\textwidth][c]{$u_{tttt}$}     &   \makebox[0.10\textwidth][c]{$u_{ttt}$}   &   \makebox[0.10\textwidth][c]{$u_{tt}$}   & \makebox[0.10\textwidth][c]{$u_{t}$} & \makebox[0.10\textwidth][c]{$u$} & \makebox[0.08\textwidth][c]{AICc}\\ \midrule
        correct     & 1.0 &  & 8.000 &  & 16.000 & -      \\
        SINDy       & 1.0  & 0.066   & 7.950 & 0.264 & 15.843  & -587.8   \\
        \textbf{LES-SINDy} & \textbf{1.0}  &  & \textbf{8.000} &  & \textbf{16.000} & \textbf{-2559.3}     \\ \bottomrule
    \end{tabular}
\end{table}

Table \ref{tab:high_order} illustrates the performance of LES-SINDy and SINDy when applied to a high-order ODE. While exploring the dynamics of high-order systems through measurements is valuable, SINDy struggles with accurately approximating high-order derivative terms, especially in the presence of noise typical of real-world sensor measurements. This limitation leads to inaccuracies in the identified models, as evidenced by the spurious terms and incorrect coefficients in the SINDy results. In contrast, LES-SINDy leverages the Laplace domain to approximate high-order derivatives more accurately, which results in a model that closely matches the true dynamics without introducing erroneous terms. Specifically, LES-SINDy accurately captures the correct coefficients for the high-order derivatives. Quantitatively, LES-SINDy improves from $-587.8$ to $-2559.3$, which demonstrates its superior performance in modeling high-order systems.

While LES-SINDy generally provides superior results in approximating \eqref{eq:homo_fourth}, we notice occasional failures during the initial tuning phase, which merit further investigation. From Fig. \ref{fig:errorbar_fourth} (right), one caveat in employing LES-SINDy is the careful selection of $s$ in the Laplace domain, which is crucial for successful model identification. When the frequencies $s$ are chosen too closely, the Laplace transformation yields similar outputs for different candidate functions, which leads to a nearly singular Laplace-enhanced matrix. With a bad construction in the Lanpalce-enhanced library, LES-SINDy is difficult to distinguish between candidate terms, which potentially results in the failure to identify the correct dynamical system. On the other hand, if the values of $s$ are too large and the corresponding time-domain function values are small, their Laplace transforms approach zero. This results in a sparse Laplace-enhanced library and potential inaccuracies in system identification. Therefore, careful selection of $s$ is essential to ensure a well-conditioned Laplace-enhanced library, and this should be factored into hyperparameter tuning.

\subsection{ODEs with discontinuous inputs}\label{subsec:exp_discontinuous}

In real-world systems, inputs often exhibit discontinuities, such as sudden jumps or impulses. These discontinuities can be modeled using step functions and delta functions, which represent instantaneous changes in the system's state or external forces acting on the system. Accurately identifying the governing equations of such systems is challenging due to the abrupt changes in dynamics, which can lead to complexities in both the modeling process and the numerical solution of these equations. We focus on two representative examples:

\begin{figure*}[!htbp]
\centering
\subfigure[ODEs with delta and step inputs.]{
\centering
\includegraphics[height=0.25\columnwidth]{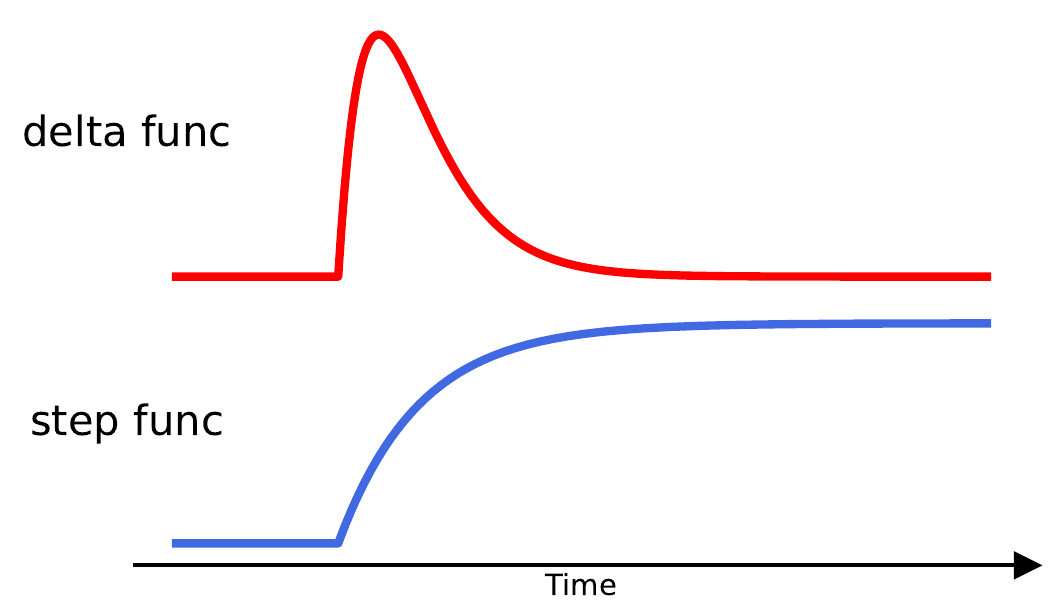}\label{fig:discontinuity}
}\ \ \ \ \ \ \ \ \ \ \ \ 
\subfigure[ODEs with trigonometric and hyperbolic functions.]
{
\centering
\includegraphics[height=0.25\columnwidth]{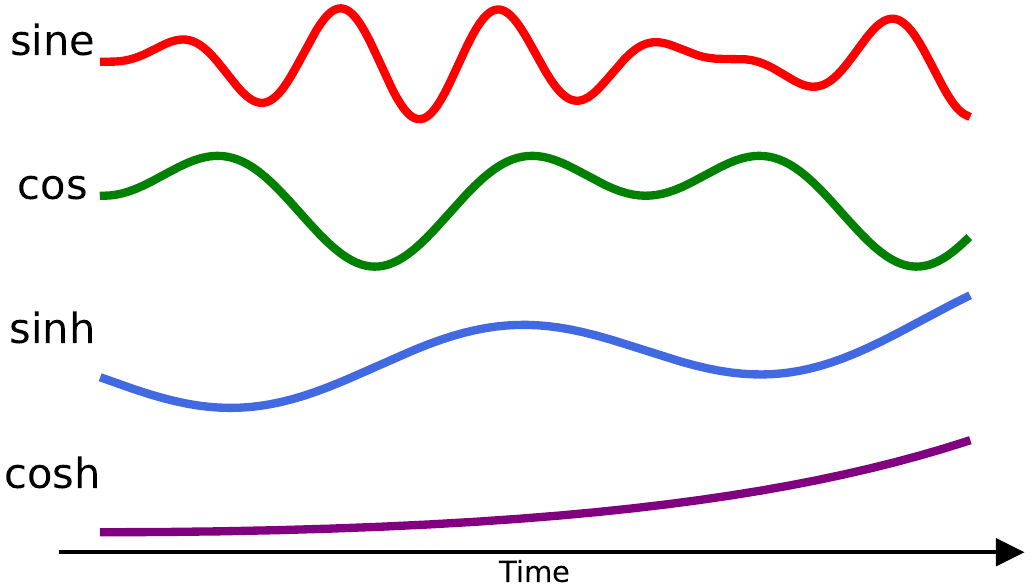}\label{fig:hyperbolic}
}

\subfigure[log RMSE error mappings.]
{
\centering
\includegraphics[height=0.40\columnwidth]{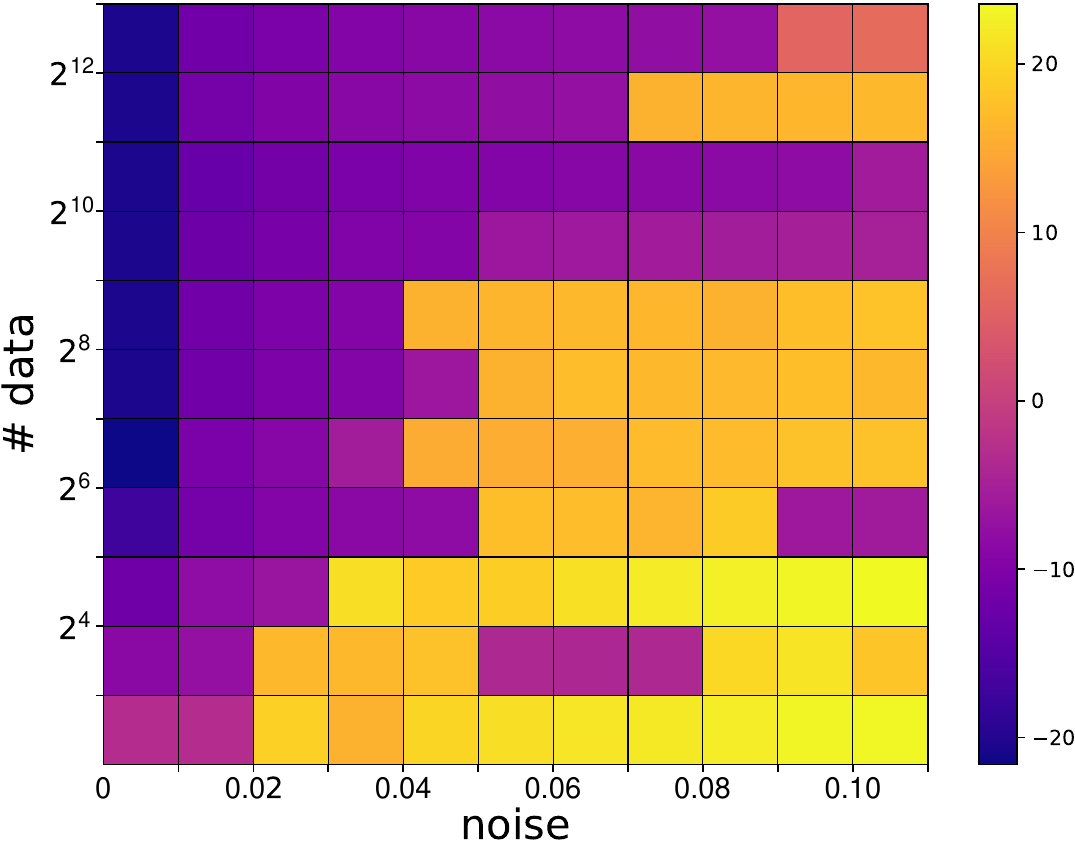}\label{fig:error_trig}\ \  
\includegraphics[height=0.40\columnwidth]{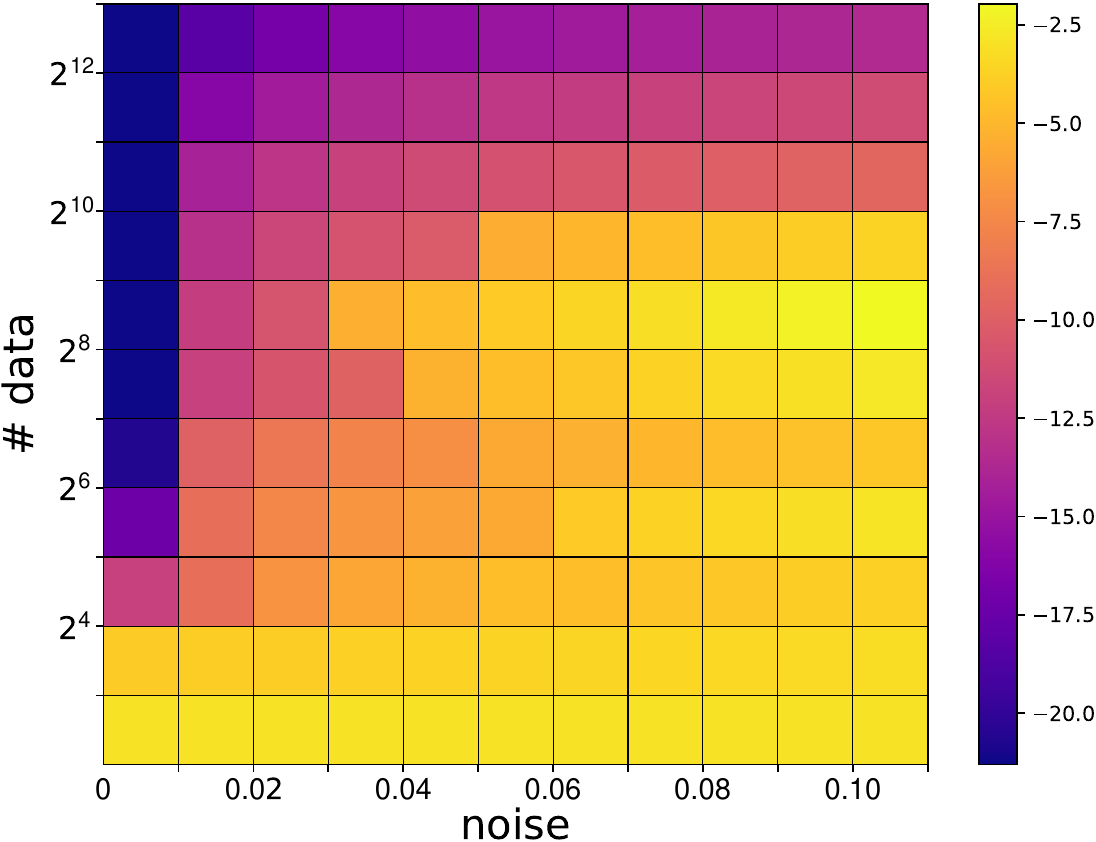}
}
% \vspace{-0.15 in}
\caption{\textbf{Top Left:} Simulation of the ODEs with Discontinuous Inputs. The x-axis denotes the time, and the y-axis represents the magnitude of the dynamical system response. \textbf{Top Right:} Simulation of the ODEs with Trigonometric (sine and cosine) and Hyperbolic (sinh and cosh) Functions. \textbf{Bottom:} log RMSE Error Mapping: Noise Scale and \# Measurements Effects (Tested in ODEs with Trigonometric Inputs) in Logarithmic Scale. The x-axis denotes the noise scale of the measurements, the y-axis represents the number of measurements used for LES-SINDy, and the colors represent the log RMSE errors between the predicted results and the ground truth. Magnitudes of noise scale increase uniformly from $0.0$ to $0.11$, and the number of measurements increases exponentially from $2^3$ to $2^{13}$. The color depth in heat mapping indicates the change in prediction accuracy as the amount of measurement noise and the amount of measurement increase. The left figure denotes the log RMSE error mapping for an ODE with a sine input, and the right one represents the results of an ODE with a cosine input.}\vspace{-0.0 in}
\end{figure*}

\textbf{ODEs with Delta function input} (shown in the upper part of Fig. \ref{fig:discontinuity}) considers a second-order linear ODE with a delta function input $\delta(\cdot)$, which models an impulsive force acting on the system. The governing equation is given by:
\begin{equation}\label{eq:ode_delta_func}
    {u}_{tt} + 2 \zeta \omega {u}_{t} + \omega^2 u = F_0 \delta(t - t_0), 
\end{equation}
where $u$ here represents the system's response, $\zeta=2.0$ is the damping ratio, $\omega=2.0$ is the natural frequency of the system, $F_0=1.0$ is the magnitude of the impulse, and $\delta(t - t_0)$ is the Dirac delta function centered at time $t_0$. This ODE is often used to model systems subjected to a sudden impact or shock, such as a mechanical system experiencing a collision or an electrical circuit subjected to a pulse.

\textbf{ODEs with Step function input} (see the bottom part of Fig. \ref{fig:discontinuity}) involves a linear ODE with a step function input $H(\cdot)$, which models a sudden and sustained change in the system's driving force. The governing equation is:
\begin{equation}
    {u}_{t} + \alpha u = F_0 H(t - t_0),
\end{equation}
where $\alpha=2.0$ is a positive constant representing the system's rate of decay or growth, $F_0=1.0$ is the magnitude of the step, and $H(t - t_0)$ is the Heaviside step function, which represents a sudden jump in the input at time $t_0$. This type of ODE is commonly used to model processes such as charging and discharging in electrical circuits or the response of mechanical systems to a sudden application of force.

We transform the delta and step functions from the time domain into the Laplace domain and incorporate them in the Laplace-enhanced library with respect to different complex frequencies $s$. The transformations are expressed as:
\begin{equation}
    \begin{split}
        \mathcal{L}\left\{\delta(t - t_0)\right\}(s) & = \int_0^\infty e^{-st}\delta(t - t_0) dt = e^{-t_0s}\\
        \mathcal{L}\left\{H(t - t_0)\right\}(s) & = \int_0^\infty e^{-st}H(t - t_0) dt = \frac{e^{-t_0s}}{s}.
    \end{split}
\end{equation}
The applications of both SINDy and LES-SINDy to discover such ODEs are summarized in Table \ref{tab:delta_func}. While SINDy effectively approximates the continuous terms, it struggles to capture the discontinuities due to its limitations in approximating discontinuous terms in the time domain. In contrast, by transforming these terms into the Laplace domain, LES-SINDy identifies the corresponding ODEs with high accuracy.

\begin{table}[!htbp]
\centering\caption{Results of LES-SINDy applied to ODEs with delta and step inputs. While SINDy fails to accurately identify the system due to inadequate approximation of discontinuities, LES-SINDy successfully identifies the correct system.}\label{tab:delta_func}
\begin{tabular}{l|ccccccc}
\toprule
\makebox[0.10\textwidth][l]{results}     &   \makebox[0.05\textwidth][c]{$u_{tt}$}     &   \makebox[0.10\textwidth][c]{$u_t$}   &   \makebox[0.10\textwidth][c]{$u$}   & \makebox[0.10\textwidth][c]{$H(t-t_0)$} & \makebox[0.10\textwidth][c]{$\delta(t-t_0)$}  & \makebox[0.08\textwidth][c]{AICc}\\ \midrule
truth     & 1.0     & 4.000   & 4.000  & & 1.000 & -   \\
SINDy       & 1.0     & 3.999   & 3.999 & & -99.980 & 126.8    \\
\textbf{LES-SINDy}   & \textbf{1.0}     & \textbf{4.000}   & \textbf{4.000} & & \textbf{1.000}   & \textbf{-1921.4}   \\ \midrule
truth     & & 1.0     & 2.000   & 1.000  & &   -    \\
SINDy       & & 1.0     & 2.000   & -1.000 & 0.505  &  423.2  \\
\textbf{LES-SINDy} &  & \textbf{1.0}     & \textbf{2.001}   & \textbf{1.000} & &  \textbf{-1468.2}    \\  \bottomrule
\end{tabular}
\end{table}

To evaluate LES-SINDy’s robustness to noise, we further introduced varying levels of noise into the time-series measurements, from clean measurements to 10\% noise. The results, as shown in Table \ref{tab:step_noise}, reveal the significant resilience of LES-SINDy to noisy measurements.

\begin{table}[!htbp]
\centering\caption{Results of LES-SINDy applied to ODE with step and delta inputs across varying levels of noise.}\label{tab:step_noise}
\begin{tabular}{l|ccccccccc}
\toprule
\makebox[0.10\textwidth][l]{results}     &   \makebox[0.08\textwidth][c]{$u_{tt}$}  &   \makebox[0.10\textwidth][c]{$u_t$}    &   \makebox[0.10\textwidth][c]{$u$}   & \makebox[0.10\textwidth][c]{$\delta(t-t_0)$} & \makebox[0.08\textwidth][c]{AICc}\\  \midrule
truth & 1.0 & 4.000 & 4.000 & 1.000 & -       \\ \midrule
clean   & 1.0 & 4.000 & 4.000 & 1.000 & -1921.4 \\
2\% noise     & 1.0 & 3.996 & 4.011 & 0.989 & -1065.1 \\
4\% noise     & 1.0 & 3.997 & 4.008 & 0.990 & -1129.3 \\
6\% noise     & 1.0 & 4.005 & 3.984 & 1.017 & -996.4 \\
8\% noise     & 1.0 & 4.006 & 3.972 & 1.037 & -901.8 \\
10\% noise    & 1.0 & 3.916 & 4.232 & 0.851 & -460.25 \\ \bottomrule
\end{tabular}
% \end{table}

% \begin{table}[!htbp]
% \centering\caption{Results of LES-SINDy applied to ODE with a delta input with varying levels of noise.}\label{tab:delta_noise}
\centering\vspace{0.05 in}
\begin{tabular}{l|ccccccc}
\toprule
\makebox[0.10\textwidth][l]{results}     &   \makebox[0.08\textwidth][c]{$u_t$}   &   \makebox[0.10\textwidth][c]{$u$}   & \makebox[0.10\textwidth][c]{$H(t-t_0)$} & \makebox[0.08\textwidth][c]{AICc}\\ \midrule
truth     & 1.0     & 2.000   & 1.000     &  - \\ \midrule
clean       & 1.0     & 2.001   & 1.000     &  -1468.2\\
1\% noise         & 1.0     & 2.002   & 1.001     &  -1186.7 \\
3\% noise         & 1.0     & 1.646   & 0.890     &  -422.1\\
5\%  noise        & 1.0     & 1.232   & 0.761     &  -225.2 \\ \bottomrule
\end{tabular}
\end{table}

% It should be noted that LES-SINDy assumes prior knowledge of discontinuity locations. In cases where this information is unavailable, the method would require the library to be augmented with numerous additional terms to detect discontinuities throughout the entire time domain. The development of an automated strategy for identifying discontinuities remains an open problem, which we intend to explore in future work.

\subsection{ODEs with trigonometric and hyperbolic functions}\label{subsec:exp_tri_and_hyper}

Trigonometric and hyperbolic functions frequently appear in contexts where oscillatory behavior, wave phenomena, or exponential growth and decay are present. These functions introduce complexity into the system's dynamics due to their periodic or rapidly varying nature. Accurately solving and analyzing such ODEs is essential for understanding systems ranging from simple harmonic oscillators to complex wave propagation and signal processing. To illustrate the application of ODEs with trigonometric and hyperbolic functions, we consider the following dynamical systems.

\textbf{ODEs with a sine function} (the top panel of Fig. \ref{fig:hyperbolic}) represents periodic oscillations. It frequently appears in models of harmonic motion and wave-like phenomena. For instance, the ODE
\begin{equation}\label{eq:sine_func}
    {u}_{tt} + 15u = 2\sin(3t)
\end{equation}
describes a system where $u$, the system response, is influenced by a sinusoidal external force characterized by $\sin(3t)$.

\textbf{ODEs with a cosine function} (the upper middle panel of Fig. \ref{fig:hyperbolic}) also models periodic behavior and is widely used in the study of vibrations and oscillations. In the equation
\begin{equation}\label{eq:cosine_func}
    {u}_{tt} + 4u = \cos(t),
\end{equation}
the cosine function $\cos(t)$ serves as the driving force, which represents mechanical vibrations or alternating current circuits. The response $u$ is subject to a periodic input.

\textbf{ODEs with a hyperbolic sine (sinh) function} (the lower middle panel of Fig. \ref{fig:hyperbolic}) is commonly associated with exponential growth and decay, where the phenomena often observed in thermal processes and certain types of wave propagation. The ODE
\begin{equation}\label{eq:sinh_func}
    {u}_{tt} + 4u = \sinh(2t)
\end{equation}
captures the dynamics of a system where the hyperbolic sine function $\sinh(2t)$ induces exponential behavior, which is relevant in contexts like signal amplification or heat distribution in hyperbolic geometries.

\textbf{ODEs with a hyperbolic cosine (cosh) function} (the bottom panel of Fig. \ref{fig:hyperbolic}) is distinct in its ability to describe stable exponential growth, which remains positive across its entire domain. The equation
\begin{equation}\label{eq:cosh_func}
    {u}_{tt} - 4u = \cosh(2t)
\end{equation}
illustrates how $\cosh(2t)$ can drive the evolution of a system, which can model scenarios where the response $u$ is governed by stable, exponential dynamics, such as the growth of instabilities in mechanical systems or certain relativistic effects.

% The analysis of these ODEs provides valuable insight into the behavior of systems influenced by oscillatory and exponential factors. Such equations serve as a foundation for understanding more complex dynamical systems and highlight the significance of trigonometric and hyperbolic functions in applied mathematical modeling.

To expand the Laplace-enhanced library with trigonometric and hyperbolic functions, we now include the Laplace transforms of these functions. Specifically, the resulting Laplace transforms for sine, cosine, sinh, and cosh functions concerning $s$ are given as follows:
\begin{equation}
    \begin{split}
        \mathcal{L}\left\{\sin(\omega t)\right\}(s) & = \frac{\omega}{s^2 + \omega^2},\\
        \mathcal{L}\left\{\cos(\omega t)\right\}(s) & = \frac{s}{s^2 + \omega^2},\\
        \mathcal{L}\left\{\sinh(\omega t)\right\}(s) & = \frac{\omega}{s^2 - \omega^2},\\
        \mathcal{L}\left\{\cosh(\omega t)\right\}(s) & = \frac{s}{s^2 - \omega^2},
    \end{split}
\end{equation}
where $\omega$ here is the angular frequency.

Table \ref{tab:func_tri_hyper} provides a summary of the correct ODEs and those identified by LES-SINDy, which illustrates its capability to handle these tasks effectively. It is noteworthy that while Fourier-enhanced SINDy can accurately approximate ODEs with sine and cosine inputs, it may fail with functions exhibiting unbounded growth, such as sinh and cosh functions in \eqref{eq:sinh_func}-\eqref{eq:cosh_func}. Intuitively, the Fourier transform assumes that candidate functions are absolutely integrable over the time domain, which presents challenges when approximating these functions. In contrast, the Laplace transform, due to its exponentially decaying weighting factor, can handle these functions by selecting $s$ with a sufficiently large real part, as demonstrated in Table \ref{tab:func_tri_hyper}. 

\begin{table}[!htbp]
\centering\caption{Results of LES-SINDy applied to ODEs with trigonometric and hyperbolic functions.}\label{tab:func_tri_hyper}
    \begin{tabular}{l|ccccccccc}
\toprule
\makebox[0.10\textwidth][l]{results}     &   \makebox[0.05\textwidth][c]{$u_{tt}$}      &   \makebox[0.08\textwidth][c]{$u$}   &   \makebox[0.08\textwidth][c]{$\sin(3t)$}   & \makebox[0.08\textwidth][c]{$\cos(t)$} & \makebox[0.08\textwidth][c]{$\sinh(2t)$} & \makebox[0.08\textwidth][c]{$\cosh(2t)$}  & \makebox[0.08\textwidth][c]{AICc} \\ \midrule
truth     & 1.0    & 15.000   & 2.000  & & & & -    \\
LES-SINDy   &  1.0   & 15.000   & 2.000  & & & & -2456.6    \\ \midrule
truth     & 1.0     & 4.000   & & 1.000 & & & -     \\
LES-SINDy   & 1.0     & 4.000   & & 1.000 & & & -1891.3     \\ \midrule
truth     & 1.0     & 4.000   & & & 1.000 & & -     \\
LES-SINDy   & 1.0     & 4.000   & & & 1.000 & & -1050.5     \\ \midrule
truth     & 1.0     & -4.000   & & & & 1.000 & -     \\
LES-SINDy   & 1.0     & -4.000   & & & & 1.000 & -2369.4     \\\bottomrule
\end{tabular}
\end{table}

We conducted additional tests on ODEs with sine and cosine inputs to examine the impact of measurement quantity and noise level on the performance of LES-SINDy. As presented in Fig. \ref{fig:error_trig}, we varied both the number of measurements and the noise level to assess whether LES-SINDy could consistently identify the correct ODEs. The measurements were uniformly sampled over the time domain. The results indicate that when the noise level exceeds 10\%, LES-SINDy often fails to identify the correct ODEs. This issue can be alleviated by increasing the number of measurements, but as noise increases linearly, the required sample size grows roughly exponentially. Therefore, both the quality and quantity of measurements collected from sensors must be carefully balanced to ensure the successful identification of the ODEs.

% Similarly to our assumption about discontinuity locations, we here assume the angular frequency is known when discovering ODEs involving trigonometric and hyperbolic functions. Without this assumption, the library must be augmented with additional terms to detect the angular frequency. We will explore the automatic detection strategy of angular frequency in future work.

\subsection{Nonlinear ODE systems}

Nonlinear ODE systems are pervasive in many areas of science and engineering, where they are used to model complex, interacting phenomena. These systems are characterized by their sensitivity to initial conditions, the presence of feedback loops, and often, chaotic behavior. Accurately identifying the governing equations of such systems is crucial for understanding and predicting their long-term behavior, which has profound implications in fields ranging from meteorology to ecology. We evaluate the performance of LES-SINDy on two iconic nonlinear ODE systems:

\begin{figure*}[!ht]
\centering
\subfigure[The Lorenz system.]{
\centering
\includegraphics[height=0.17\columnwidth]{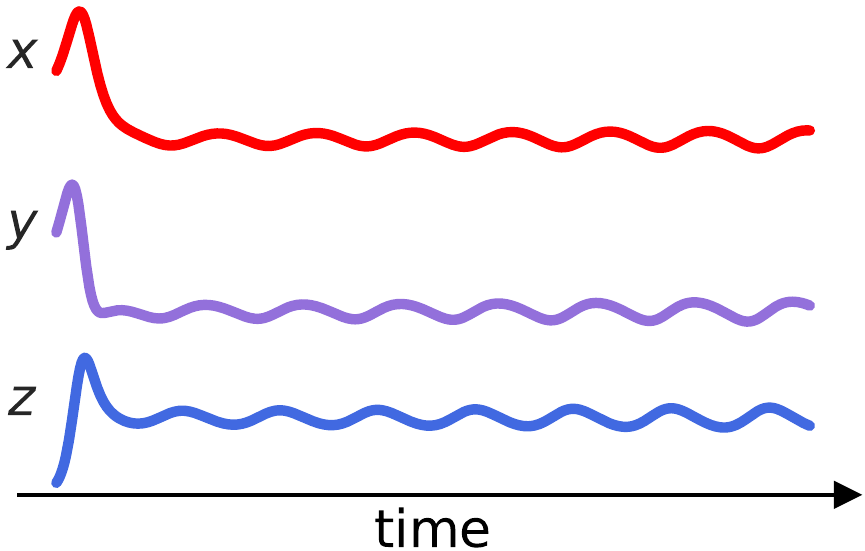}\label{fig:lorenz2}
\hspace{-0.075 in}
\includegraphics[width=.3\columnwidth]{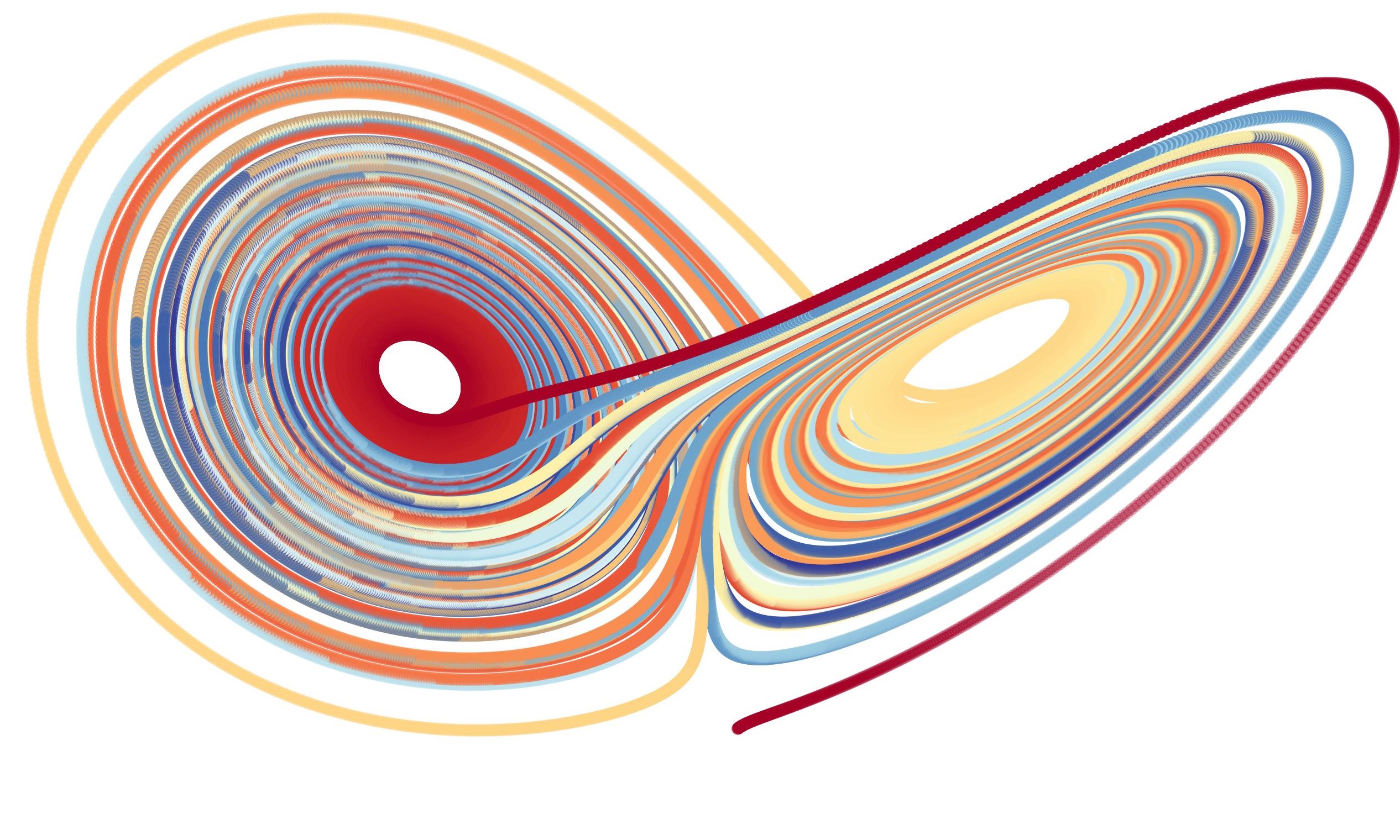}\label{fig:lorenz}
}\ \ \ \ 
\subfigure[The Lotka–Volterra model.]
{
\centering
\includegraphics[height=0.17\columnwidth]{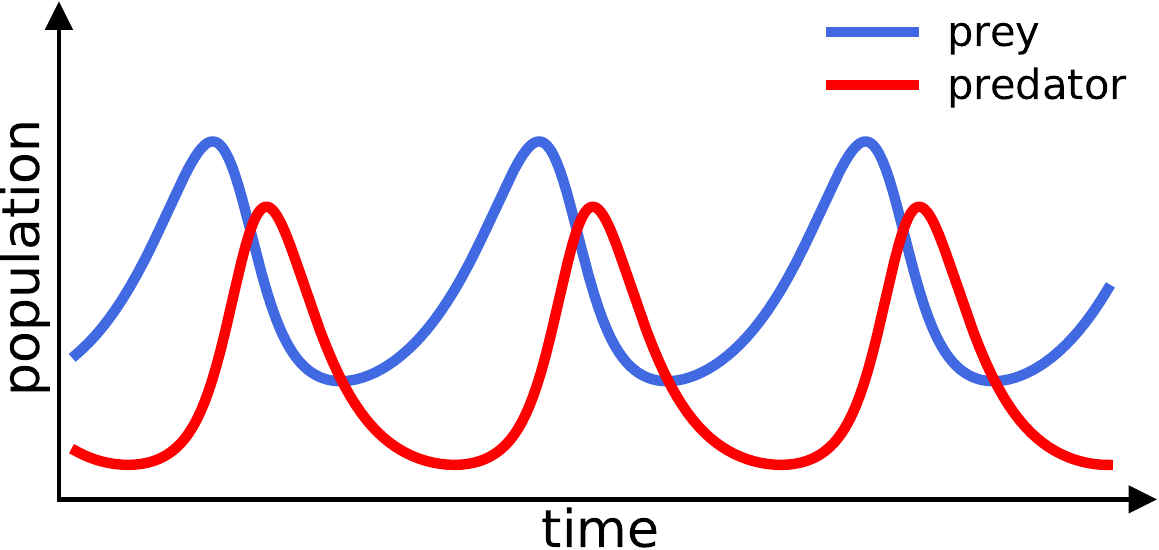}\label{fig:lv}
}
% \vspace{-0.15 in}
\caption{Simulation of Nonlinear ODE systems. \textbf{Left:} Simulation of the Lorenz system's state variables $u=\left[x, y, z\right]^\intercal$. The x-axis denotes the time, and the y-axis represents the magnitude of the dynamical system. \textbf{Middle:} The corresponding trajectory of the Lorenz system in the three-dimensional phase space. \textbf{Right:} Simulation of the Lotka–Volterra model. The x-axis denotes the time, and the y-axis represents the population dynamics of prey and predator.}
\label{fig:ode_system}\vspace{-0.0 in}
\end{figure*}

\textbf{The Lorenz system} (Fig. \ref{fig:lorenz}) is one of the most famous examples of a chaotic system, originally derived from simplified equations of atmospheric convection. It is governed by a set of three coupled nonlinear ODEs:

\begin{equation}\label{eq:lorenz}
    \begin{aligned}
   x_t &= \sigma (y - x), \\
   y_t &= x (\rho - z) - y, \\
   z_t &= x y - \beta z,
   \end{aligned}
\end{equation}
where $x$, $y$, and $z$ ($\boldsymbol{u}(t)\in \left[x(t), y(t), z(t)\right]^\intercal$ and normally we omit $t$ for simplicity) represent the system's state variables, and \(\sigma=10.0\) (the Prandtl number), \(\rho=28.0\) (the Rayleigh number), and \(\beta=8/3\) (the aspect ratio) are parameters derived from the physical properties of the system. These parameters are not arbitrary but are derived from physical properties in the system that Lorenz was modeling, which is a simplified model of atmospheric convection. With a given set of model parameters, the Lorenz system is well-known for its chaotic solutions, which are highly sensitive to initial conditions. Testing LES-SINDy on this system allows us to evaluate its ability to identify the underlying dynamics in the presence of chaotic behavior, where small errors in the model can lead to vastly different predictions over time.

To evaluate the performance of LES-SINDy, we applied it to time-series measurements generated by the Lorenz system. The initial model learned by LES-SINDy was expressed as:
\begin{equation}\label{eq:lorenz_learn1}
    \begin{aligned}
       &1.000z+0.375z_t-0.375xy=0, \\
       &-23.974x-3.015y+0.401x_t+1.000y_t+1.000xz=0, \\
       &-3.180x+3.179y+2.667z-0.318x_t+1.000z_t-1.000xy=0,
   \end{aligned}
\end{equation}

This learned model captures the key interactions between the variables but requires further refinement to closely match the canonical form of the Lorenz equations. By applying a suitable variable transformation, the learned model was refined to the following system:
\begin{equation}\label{eq:lorenz_learn2}
    \begin{aligned}
        & x_t + 10.007x - 10.006 y = 0,\\
	& y_t - 27.989x + 0.998 y + 1.000xz = 0,\\
        & z_t - 1.000xy + 2.667 z = 0,\\
   \end{aligned}
\end{equation}
which is remarkably close to the original Lorenz system. The identified equations exhibit only minor deviations in the coefficients from the true parameters. The close match between the learned and true models demonstrates the potential of LES-SINDy to be a powerful tool for discovering the underlying dynamics of chaotic systems.

\textbf{The Lotka-Volterra model} (Fig. \ref{fig:lv}), which is also referred to as the predator-prey model, consists of a pair of first-order nonlinear differential equations that describe the interaction between two species: a predator and its prey. The Lotka-Volterra model is mathematically represented by a system of nonlinear ODEs:
\begin{equation}\label{eq:lotka_volterra}
   \begin{aligned}
   x_t &= \alpha x - \beta xy, \\
   y_t &= \delta xy - \gamma y,
   \end{aligned}
\end{equation}
where $x$ represents the prey population, and $y$ represents the predator population ($\boldsymbol{u}(t)\in \left[x(t), y(t)\right]^\intercal$). The model depends on four parameters: the prey growth rate ($\alpha=1.0$), which represents the exponential growth of the prey population when predators are absent; the predation rate coefficient ($\beta=1.0$) governs the dynamics between predator and prey populations and reflects the efficiency with which predators consume prey; the predator mortality rate ($\delta=1.0$) indicates how quickly the predator population decreases when no prey is available; and the predator reproduction rate ($\gamma=1.0$) describes the predator population growth as a result of prey consumption.

\begin{figure}[!htbp]
    \centering
    {\includegraphics[width=1.0\linewidth]{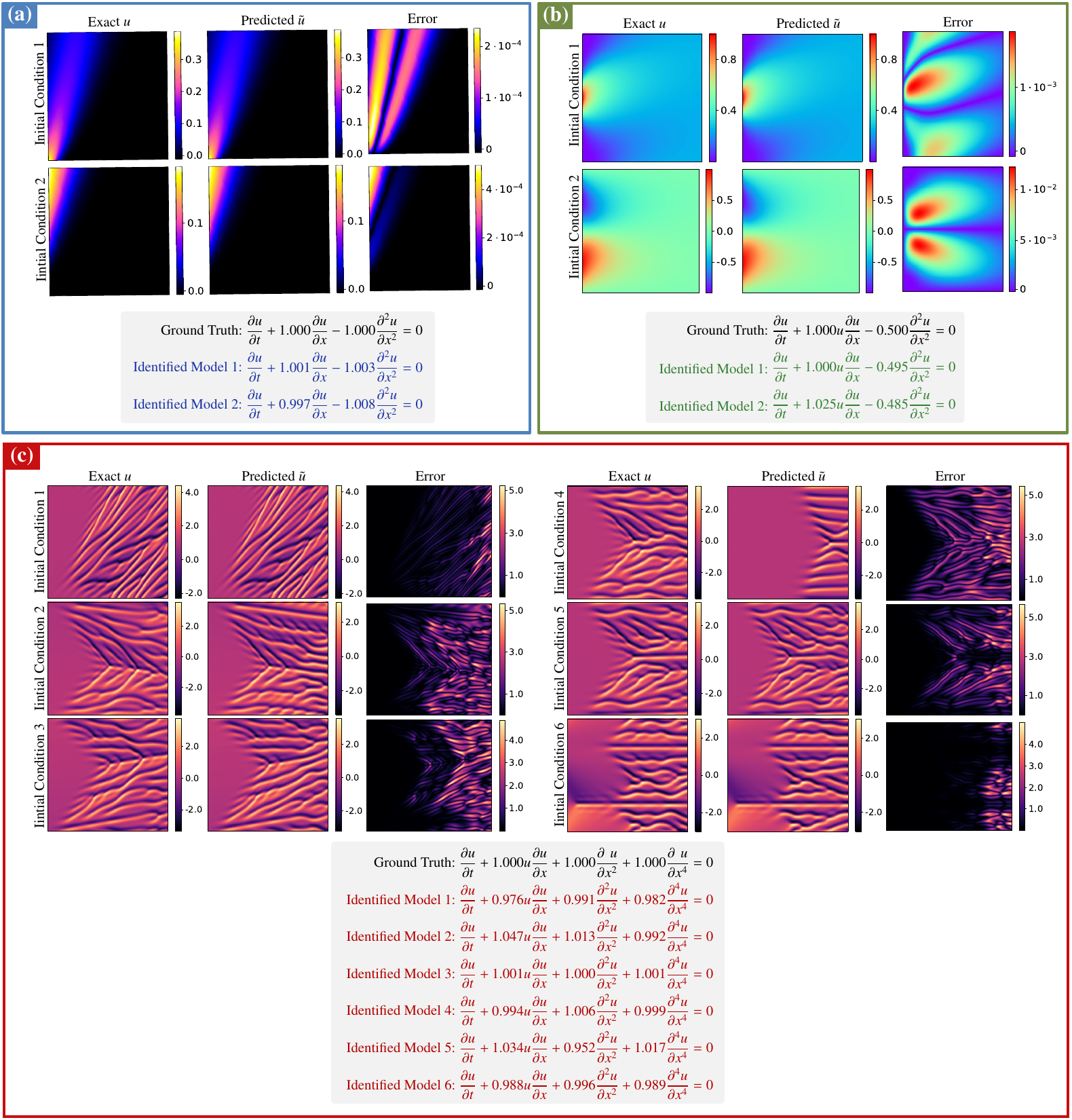}}
    \caption{Identified models compared with the ground truth for canonical PDEs, tested with different initial conditions.  \textbf{Top Left:} the convection-diffusion equation. \textbf{Top Right:} the Burgers equation. \textbf{Bottom:} Kuramoto-Sivashinsky equation.}
    \label{fig:pde_results}
\end{figure}

The performance of LES-SINDy was further evaluated through its application to time-series measurements from the Lotka-Volterra system. The identified model:
\begin{equation}\label{eq:lv_learn1}
   \begin{aligned}
   -1.004x-5.071y+1.000x_t-5.071y_t+6.072xy=0, \\
   0.195x-1.000y-0.194x_t+1.000y_t-19194xy=0,
   \end{aligned}
\end{equation}
succeeds in depicting predator-prey interactions. Nevertheless, the model's coefficients and structure needed adjustment to more accurately reflect the canonical Lotka-Volterra equations. Upon applying an appropriate variable transformation, the model was simplified to a more manageable form:
\begin{equation}\label{eq:lv_learn2}
   \begin{aligned}
   x_t - 1.009 x + 1.004 xy=0, \\
   y_t + 1.000 y - 1.000 xy=0.
   \end{aligned}
\end{equation}
The refined model is closely related to the original Lotka-Volterra equations, with only slight variations in the coefficients. These results confirm that the identified equations accurately represent the predator-prey dynamics, even in the presence of intricate, nonlinear interactions.
% \textcolor{red}{AICc eliminate redundant variables, only test one ODE}

\subsection{Partial differential equations}\label{subsec:pde}

PDEs play a fundamental role in modeling a wide range of physical phenomena that involve spatial and temporal variations. These include processes such as fluid flow, heat conduction, wave propagation, and pattern formation. PDEs are inherently more complex than ODEs because they involve derivatives with respect to both time and space. Accurately identifying the governing PDEs from measurements is crucial for predicting and controlling the behavior of these systems in various scientific and engineering applications. We apply LES-SINDy to three representative PDEs that are widely studied in applied mathematics and physics.

\textbf{The Convection-Diffusion equation} (Fig. \ref{fig:pde_results}\textcolor{darkred}{(a)}) is a fundamental PDE that models the transport of a scalar quantity (such as heat, mass, or pollutants) in a medium. It combines the effects of advection (transport due to a velocity field) and diffusion (spreading due to concentration gradients) and is given by:
\begin{equation}\label{eq:advection}
   \frac{\partial u}{\partial t} + c \frac{\partial u}{\partial x} - D \frac{\partial^2 u}{\partial x^2} = 0,
\end{equation}
where $u(x,t)$ represents the quantity being transported, $c=1.0$ is the advection speed, and $D=1.0$ is the diffusion coefficient. The convection-diffusion equation is essential in fields such as fluid dynamics, environmental science, and chemical engineering. Testing LES-SINDy on this equation allows us to evaluate its ability to identify PDEs that involve both first-order and second-order spatial derivatives.

\begin{table}[!htbp]
\centering\caption{Identification of the convection-diffusion equation across different novel levels using LES-SINDy. It performs accurately under noise conditions ranging from 0\% to 30\%.}\label{tab:advection_result}
      \begin{tabular}{l|ccccccc}
		\toprule
            \makebox[0.10\textwidth][l]{results}     &   \makebox[0.08\textwidth][c]{$u_{t}$}  &   \makebox[0.10\textwidth][c]{$u_x$}    &   \makebox[0.10\textwidth][c]{$u_{xx}$}  & \makebox[0.08\textwidth][c]{AICc}   \\  \midrule
            truth & 1.0 & 1.000 & -1.000  & - \\ \midrule
		clean & 1.0 & 1.000 & -1.000  & -1931.5 \\
		2\% noise   & 1.0 & 1.000 & -1.000  & -1461.7 \\
		5\% noise   & 1.0 & 1.000 & -1.000  & -1226.6 \\
		10\% noise   & 1.0 & 0.986 & -0.996  & -1028.7 \\
		20\% noise   & 1.0 & 0.969 & -0.990  & -945.4 \\ 
		30\% noise   & 1.0 & 0.951 & -0.983  & -898.6 \\ \bottomrule
		\end{tabular}
\end{table}

\textbf{Burgers equation}, shown in Fig. \ref{fig:pde_results}\textcolor{darkred}{(b)}, is a nonlinear PDE that serves as a simplified model for various physical processes, including fluid dynamics, traffic flow, and gas dynamics. It is often used as a test case for numerical methods due to its balance of nonlinearity and diffusion. The equation is expressed as:
\begin{equation}\label{eq:burgers}
    \frac{\partial u}{\partial t} + u \frac{\partial u}{\partial x} - \nu \frac{\partial^2 u}{\partial x^2} = 0,
\end{equation}
where $u(x,t)$ is the velocity field, and $\nu=0.50$ is the viscosity coefficient. Burgers equation is notable for its ability to develop shock waves (discontinuities in the solution), which makes it a challenging test for identification. LES-SINDy’s performance on Burgers equation will demonstrate its effectiveness in capturing nonlinear dynamics and shock formation in PDEs.

\begin{table}[!htbp]
\centering\caption{Identification of the Burgers equation across different novel levels using LES-SINDy. It can perform accurately under noise conditions ranging from 0\% to 20\%.}\label{tab:burgers_result}
    \begin{tabular}{l|ccccccc}
		\toprule
            \makebox[0.10\textwidth][l]{results}     &   \makebox[0.08\textwidth][c]{$u_{t}$}  &   \makebox[0.10\textwidth][c]{$u_{xx}$}    &   \makebox[0.10\textwidth][c]{$uu_{x}$}  & \makebox[0.08\textwidth][c]{AICc}   \\  \midrule
            truth & 1.0 & -0.500 & 1.000 & - \\ \midrule
		clean   & 1.0 & -0.500 & 1.001 & -1360.1  \\
		2\% noise     & 1.0 & -0.498 & 0.999 & -674.0  \\
		5\% noise     & 1.0 & -0.492 & 0.996 & 718.9  \\
		10\% noise    & 1.0 & -0.503 & 1.058 & 993.9  \\
		20\% noise    & 1.0 & -0.504 & 1.070 & 1189.5  \\ \bottomrule
		\end{tabular}
\end{table}

\textbf{The Kuramoto-Sivashinsky (KS) equation} in Fig. \ref{fig:pde_results}\textcolor{darkred}{(c)} is a fourth-order nonlinear PDE that arises in the study of instabilities in fluid flows, flame fronts, and other dissipative systems. It is known for its chaotic behavior and pattern formation, which encapsulates a rich array of phenomena including pattern formation, chaos, and turbulence within a relatively simple mathematical framework. Its complex spatiotemporal behaviors enable the exploration of dissipative systems far from equilibrium, which also makes it a challenging model for both analysis and identification. The KS equation is given by:
\begin{equation}\label{eq:ks}
       \frac{\partial u}{\partial t} + u \frac{\partial u}{\partial x} + \frac{\partial^2 u}{\partial x^2} + \frac{\partial^4 u}{\partial x^4} = 0,
\end{equation}
where $u(x,t)$ describes the amplitude of the instability. The equation contains both second-order and fourth-order spatial derivatives, as well as nonlinear advection, which leads to complex spatiotemporal dynamics. Applying LES-SINDy to the KS equation can further demonstrate its ability to handle higher-order derivatives and chaotic behavior in spatially extended systems.

In the context of PDEs, approximating spatial derivatives of higher order poses a considerable challenge, as these derivatives are highly susceptible to noise and boundary-related issues. These difficulties often cause substantial inaccuracies in identifying the correct dynamical systems. To mitigate this challenge when identifying PDEs, we have to modify the construction of the Laplace-enhanced library to suit the nature of PDEs. Instead of following previous approaches that apply the Laplace transform over the time domain, it is more effective here to use the transform over the spatial domain. This adjustment improves the approximation of high-order derivative terms, which minimizes the instability and errors typically induced by noise, thereby ensuring a more accurate and robust identification of the system.

The results of LES-SINDy and baseline models are presented in Fig. \ref{fig:pde_results} and Tables \ref{tab:advection_result}-\ref{tab:ks_result}. We first identified the convection-diffusion equation, Burgers equation, and KS equation under different initial conditions. As illustrated in Fig. \ref{fig:pde_results}, LES-SINDy consistently produces robust models across these different initial conditions. The close match between the discovered equations and the ground truths across all cases demonstrates the potential of LES-SINDy for uncovering underlying dynamics in complex systems. For the convection-diffusion equation, time-series measurements with noise levels were considered ranging from 0\% to 30\%. Although the identification becomes difficult as noise accumulates,  results recorded in Table \ref{tab:advection_result} indicate that LES-SINDy can reliably identify the convection-diffusion equation across increasing noise levels. The AICc values increase with higher noise levels, which reflects the reduced goodness-of-fit due to noise contamination. However, the method consistently identifies the correct model structure and provides accurate coefficient estimates despite the noise. These results underscore the effectiveness of LES-SINDy in accurately recovering the underlying dynamics of the convection-diffusion equation under varying noise conditions.

We subsequently identified the Burgers equation, with noise levels increasing from 0\% to 20\%. The results are shown in Table \ref{tab:burgers_result}, which further demonstrates the ability of LES-SINDy to discover complex PDEs in a noisy environment. 

\begin{table}[!htbp]
\centering\caption{Identification of the Kuramoto-Sivashinsky equation using LES-SINDy and baselines. LES-SINDy demonstrates superior accuracy and robustness across noise levels ranging from 0\% to 20\%.}\label{tab:ks_result}
  \begin{tabular}{l|l|ccccc}
\toprule
\multirow{2}{*}{noise} & \multirow{2}{*}{method} & \multicolumn{5}{c}{results}            \\
                       &                         & $u_t$ & $u_{xx}$ & $u_{xxxx}$ & $uu_x$  & \makebox[0.08\textwidth][c]{AICc}  \\ \midrule
\multicolumn{2}{c|}{ground truth}                        & 1.0   & 1.0      & 1.0        & 1.0   & -    \\ \midrule
\multirow{3}{*}{clean} & SINDy                   & 1.0 & 0.984    & 0.994      & 0.982   & 1268.6    \\
                       & {Weak-SINDy}     & 1.0 & {0.997} & {0.992}  & {0.995} & {1170.5}   \\
                       & LES-SINDy               & 1.0 & 0.990    & 0.995      & 0.988 & 1189.3   \\ \midrule
1\%                    & SINDy                   & 1.0 & 0.459    & 0.481      & 0.492  & 2079.6  \\ \midrule
\multirow{2}{*}{5\%}   & Weak-SINDy              & 1.0 & \textbf{0.992}    & 0.979      & 0.953 & 1353.1 \\
                       & \textbf{LES-SINDy}      & 1.0 & 0.990    & \textbf{0.985}      & \textbf{0.989} & \textbf{1306.9}  \\ \midrule
\multirow{2}{*}{10\%}  & Weak-SINDy              & 1.0 & 0.959    & 0.923      & 0.791  & 1834.7\\
                       & \textbf{LES-SINDy}     & 1.0 & \textbf{0.989} & \textbf{0.984}  & \textbf{0.986} & \textbf{1318.6}\\ \midrule
\multirow{2}{*}{20\%}  & Weak-SINDy              & 1.0 & 0.886    & 0.786      & 0.405  &  2316.6 \\
                       & \textbf{LES-SINDy}     & 1.0 & \textbf{0.957} & \textbf{0.973}  & \textbf{0.954} & \textbf{1445.2} \\ \bottomrule
\end{tabular}
\end{table}

Finally, the KS equation was tested with noise levels ranging from 0\% to 20\% and compared with baselines such as SINDy and Weak-SINDy. The results are summarized in Table \ref{tab:ks_result}. Under noise-free conditions, all three methods recovered the coefficients of the KS equation, with LES-SINDy and Weak-SINDy showing slightly higher accuracy. As noise increased, the performance of the methods began to diverge significantly: At 1\% noise, SINDy exhibited a marked decline in performance, with coefficient estimates dropping to approximately 0.5 for $u_{xx}$, $u_{xxxx}$, and $u u_x$. This indicates a high sensitivity to even small levels of noise, which renders SINDy less effective in scenarios where measurement quality cannot be guaranteed. At 5\% noise, both Weak-SINDy and LES-SINDy remained highly accurate. Weak-SINDy estimated the coefficients as 0.992 for $u_{xx}$, 0.979 for $u_{xxxx}$, and 0.953 for $u u_x$. LES-SINDy produced similar results with coefficients of 0.990, 0.985, and 0.989, respectively. LES-SINDy had a slightly lower AICc value of 1306.9 compared to Weak-SINDy's 1353.1, which indicates a better model fit. At 10\% noise, LES-SINDy demonstrated greater robustness compared to Weak-SINDy, which accurately identifies the model with an AICc of 1318.6, whereas Weak-SINDy yields a higher AICc of 1834.7. At 20\% noise, the advantages of LES-SINDy became more pronounced: LES-SINDy maintained accurate coefficient estimates of 0.957 for $u_{xx}$, 0.973 for $u_{xxxx}$, and 0.954 for $u u_x$, along with a lower AICc of 1445.2. In contrast, Weak-SINDy's coefficients degraded to 0.886, 0.786, and 0.405, with a significantly higher AICc of 2316.6. As demonstrated in Table \ref{tab:ks_result}, LES-SINDy could tolerate noise up to 20\%, whereas the other methods failed under the same conditions.

By extending LES-SINDy to these three representative PDEs, we demonstrate its versatility and robustness in identifying the governing equations of complex systems that exhibit both temporal and spatial dynamics. The successful identification of such PDEs with LES-SINDy further implies its potential for advancing our understanding of complex dynamical systems.

\section{Discussion}\label{sec_discuss}

In this study, we propose LES-SINDy, a SINDy-inspired framework enhanced by the Laplace transform, which addresses the challenges of high-order derivatives, discontinuities, and unbounded growth in functions. By applying the Laplace transform and integration by parts, LES-SINDy is systematically evaluated across diverse dynamical systems and compared its effectiveness to other established methods. The experimental findings indicate that LES-SINDy not only supports accurate data-driven model discovery in traditional tasks but also effectively captures complex dynamics, including non-smooth and exponentially growing functions, which are often challenging to capture with standard SINDy approaches.

We further discuss several challenges encountered in LES-SINDy, particularly when involving noisy measurements, the influence of coefficient magnitudes in dynamical systems, the selection of complex frequencies, and the impact of measurement quality and quantity. LES-SINDy consistently performs well compared to other approaches when dealing with noisy measurements. However, if system parameters vary widely in magnitude, the method may incorrectly eliminate important coefficients. The selection of the complex frequency $s$ is crucial for a successful model discovery process: its real part must be sufficiently large to minimize the influence of potential measurements taken after the sampled period. This ensures an accurate approximation of candidate functions in the Laplace domain. As noise levels increase, LES-SINDy may struggle to discover the correct system. Increasing the number of measurements can mitigate this, but the required sample size may grow exponentially with noise. Therefore, it is crucial to balance the quality and quantity of measurements for accurate model discovery.

Our work also suggests several promising directions for future investigation, some of which are not fully addressed in this paper. One area of interest is the handling of ODEs with discontinuities, trigonometric, or hyperbolic functions. We currently assume that the locations of these discontinuities are known. However, if this information is not available, the library needs to include many additional terms to detect discontinuities over the entire time domain. Similarly, we currently assume the angular frequency in the trigonometric and hyperbolic functions is known. Future research could focus on developing more comprehensive libraries or automatic strategies for detecting discontinuities and identifying angular frequencies without prior knowledge.

Another compelling direction involves the application of LES-SINDy to the data-driven model discovery of coordinates and governing equations. When measurements are taken from high-dimensional spaces, neural networks could help map these to lower-dimensional representations, thereby extracting essential features for discovering underlying dynamical systems. While existing work \cite{champion2019data} has leveraged automatic differentiation to achieve this, the corresponding methodologies in LES-SINDy are still under-explored. Future studies will aim to integrate neural networks with LES-SINDy to enhance the discovery of complex dynamical systems from high-dimensional measurements.

\section{Problem Statement}

% In this section, we focus on the problem of identifying implicit DEs and present LES-SINDy as a solution. We first generate a library of candidate functions in the time domain and then transform this library into the Laplace domain. Next, we apply sparse regression to determine the sparse coefficients to identify models. Finally, the models are evaluated using log RMSE and AICc to produce the most accurate and parsimonious dynamical systems are selected. For clarity, we provide a generalized overview of the workflow in Fig. \ref{fig:framework} and Algorithm \ref{alg:les-sindy}. 

% \subsection{Problem statement}
The problem of SINDy with implicit DEs focuses on the form of the DE governing the system's dynamics from the given measurements, which may not be explicitly known a priori. Specifically, given time series measurements of the state variables $\boldsymbol{u}(t)$ and their derivatives, the objective is to identify a sparse function $f$ that can be succinctly expressed in its general form as:
$$
f\left(t, \boldsymbol{u}(t), \frac{\mathrm{d} \boldsymbol{u}(t)}{\mathrm{d} t}, \frac{\mathrm{d}^2 \boldsymbol{u}(t)}{\mathrm{d} t^2}, \dots, \frac{\mathrm{d}^k \boldsymbol{u}(t)}{\mathrm{d} t^k} \right) = 0.
$$
Here, $f$ implicitly describes the system's dynamics by capturing the relationships between time $t$, the state variables $\boldsymbol{u}(t)\in\mathbb R^d$, and its derivatives up to the $k$-th order with respect to \(t\), with a focus on minimizing the number of nonzero terms in $f$.

The choice to define the problem using implicit DEs \cite{zhang2018robust, kaheman2020sindy} rather than the explicit form \cite{brunton2016discovering} provides a broader and more flexible framework. This configuration relies on less stringent assumptions about the system structure, which makes it more versatile and applicable to a wide range of dynamical systems. In contrast, the explicit form assumes a simpler relationship between variables, which can be inadequate for capturing the dynamics of more intricate systems. For instance, the Laguerre DE $\left(t \frac{\mathrm{d}^2 u}{\mathrm{d} t^2} + (1-t) \frac{\mathrm{d} u}{\mathrm{d} t} + u = 0\right)$ and the Van der Pol oscillator $\left(\frac{\mathrm{d}^2 u}{\mathrm{d} t^2} - \mu \left(1 - u^2\right) \frac{\mathrm{d} u}{\mathrm{d} t} + u = 0\right)$ exemplify cases where the explicit form is insufficient and highlight the need for a more general form.

\section{Methods}\label{sec_method}

In this section, we focus on the problem of identifying implicit DEs and present LES-SINDy as a solution. We first generate a library of candidate functions in the time domain and then transform this library into the Laplace domain. Next, we apply sparse regression to determine the sparse coefficients to identify models. Finally, the models are evaluated using log RMSE and AICc to produce the most accurate and parsimonious dynamical systems are selected. For clarity, we provide a generalized overview of the workflow in Fig. \ref{fig:framework} and Algorithm \ref{alg:les-sindy}.

\subsection{Library for implicit differential equations}\label{subsec:tensor_product}

To determine the function from measurements collected from sensors, we collect a time history of the state (time-series measurements) $\boldsymbol u(t)=\left[u_1(t),u_2(t),\cdots,u_d(t)\right]^\intercal\in \mathbb R^{d}$. The measurements are sampled at several time $t_1,t_2,\cdots,t_m \in\mathbb R_{\geq0}$ and are further arranged into a matrix with  $1$ (constant term) and $t$ (time term):

\begin{equation}\label{eq:measurements}
% \begin{eqnarray}\label{eq:measurements}
 \begin{bmatrix} 1 & t_1 &\bu^\intercal(t_1)\\ 1 & t_2 &\bu^\intercal(t_2) \\ \vdots & \vdots & \vdots \\ 1 & t_m &\bu^\intercal(t_m)\end{bmatrix} 
= \overset{\text{\normalsize state}}{\left.\overrightarrow{\begin{bmatrix}
1 & t_1 & u_1(t_1) & u_2(t_1) & \cdots & u_d(t_1)\\
1 & t_2 & u_1(t_2) & u_2(t_2) & \cdots & u_d(t_2)\\
\vdots & \vdots &\vdots & \vdots & \ddots & \vdots \\
1 & t_m & u_1(t_m) & u_2(t_m) & \cdots & u_d(t_m)
\end{bmatrix}}\right\downarrow}\begin{rotate}{270}\hspace{-.125in}time~~\end{rotate}
\end{equation}

To identify implicit DEs from the given measurements, we construct a comprehensive library of candidate functions. With the operator $\mathcal D^k: \mathbb R^{d}\rightarrow \mathbb R^{kd+d+2}$, we generate a vector that includes constant term, time term, and state vector $\bu(t)$ along with the time derivatives up to the $k$-th order of $\bu(t)$:
$$
{\mathcal D}^k \bu(t) = \left[1\ \ t\ \ u_1(t)\ \cdots\ u_d(t)\ \ \frac{\mathrm{d} u_1(t)}{\mathrm{~d} t}\ \cdots \ \frac{\mathrm{d}^k u_1(t)}{\mathrm{~d} t^k}\ \cdots \ \frac{\mathrm{d}^k u_d(t)}{\mathrm{~d} t^k}\right]^\intercal.
$$

% \textcolor{red}{Simplify the notation, use $\mathcal D$ to represent ${\mathcal D}^k \bu^\intercal(t)$. 
% Also, be careful about the \textbf{Transpose}!!} \textcolor{red}{Whether use $\frac{du}{\mathrm{~d}t}$ or $\dot u$ or $u^{(k)}$? Keep consistent.} 
The construction of the operator $\mathcal{D}^k$ is intentionally versatile, which enables it to include a comprehensive set of candidate functions for constructing implicit DEs. In its standard form, the operator comprises constant terms, time terms, state vectors, and their high-order derivatives to represent a parsimonious dynamical system. However, depending on the specific DEs being analyzed, such as those involving discontinuous inputs (e.g., delta or step functions), trigonometric or hyperbolic functions, or spatial derivatives in PDEs, the operator is flexible and can be adapted to include additional terms that capture these features. This adaptability ensures that the operator can be tailored to accurately capture the dynamics of more complex systems. To keep the explanation straightforward, we focus on a basic configuration to provide more intuitions, with further elaboration on different configurations provided in Section \ref{sec_exp}.

Building upon the operator, we further extend the candidate function by introducing the $n$-th tensor product of the candidate functions generated by $\mathcal D$:
\begin{equation*}
	\begin{split}
		\bigotimes^n{\mathcal D} & = \underbrace{\mathcal D \otimes \mathcal D \otimes \cdots \otimes \mathcal D}_{n \text{ times}} \\
            & = \bigotimes^n\left[1\ \ t\ \ u_1(t)\ \cdots\ u_d(t)\ \ \frac{\mathrm{d} u_1(t)}{\mathrm{~d} t}\ \cdots \ \frac{\mathrm{d}^k u_1(t)}{\mathrm{~d} t^k}\ \cdots \ \frac{\mathrm{d}^k u_d(t)}{\mathrm{~d} t^k}\right]^\intercal \\
		& = \left[1\ \ t\ \ u_1(t)\ \cdots\ \left(u_n(t)\right)^n\ \cdots \ \left(\frac{\mathrm{d} u_d(t)}{\mathrm{~d} t}\right)^n\ \cdots\ \left(\frac{\mathrm{d}^k u_d(t)}{\mathrm{~d} t^k}\right)^n\right]^\intercal, \\
	\end{split}
\end{equation*}
where $n$ denotes the number of tensor product copies of the candidate function in ${\mathcal D}$, and $\bigotimes^n: \mathbb R^{kd+d+2}\rightarrow \mathbb R^{\mathbf{d}}$ with $\mathbf{d}=\frac{({kd+d+n+1})!}{n!({kd+d+1})!}$ serves as the tensor product operator, which yields a comprehensive set of candidate functions in the time domain. To clarify the output generated by $\bigotimes^n{\mathcal D}$, we provide some examples in Appendix \ref{appendix:tensor_candidate} to illustrate the construction. 

While tensor products are powerful tools for capturing complex dynamical systems across varying levels of complexity, the number of candidate functions increases polynomially as the parameters $d$ (state dimensions), $k$ (derivative orders), and $n$ (number of tensor product copies) grow. This polynomial growth makes the process of model discovery more challenging. To streamline this process, additional system knowledge can be used to exclude irrelevant candidate functions before the model discovery. Also, the parameters $d$, $k$, and $n$ can be incrementally adjusted to explore candidate functions in a sequential manner, which starts with lower-order terms. The model discovery process concludes when the error metric is sufficiently small, which indicates the successful discovery of the governing physical laws.

We denote $\mathtt{X} = \mathtt{X}(d,k,n,t):= \bigotimes^n \mathcal{D}$ to represent the vector output from the tensor product. For simplicity, the parameters $d$, $k$, $n$, and $t$ will be omitted unless necessary. The individual elements of the vector $\mathtt{X}$ are denoted as $\mathtt{X}_1, \mathtt{X}_2, \dots, \mathtt{X}_{\mathbf{d}}$, and use $\mathtt{X}_{-i}$ $(i=1,2,\dots,\mathbf{d})$ to denote the vector $\mathtt{X}$ except the $i$-th element $\mathtt{X}_i$. To identify the dynamical system, a common approach \cite{brunton2016discovering, zhang2018robust, kaheman2020sindy} involves finding a vector $\zeta_i \in \mathbb{R}^{\mathbf{d}-1}$ such that, for measurements taken at different time from $t=t_1$ to $t=t_m$, the vector aligns the left-hand side $\mathtt{X}_{-i} \zeta_i$ with the right-hand side $\mathtt{X}_i$ as closely as possible:
\begin{equation}\label{eq:solve_tls}
       \mathtt{X}_{-i} \zeta_i = \mathtt{X}_i.
\end{equation}
This process is iterated from $i = 1$ to $i = \mathbf{d}$. Subsequently, the identified models are evaluated using a specific metric, and the best-performing model is selected as the identified dynamical system.

While this approach has been developed to discover dynamical systems in a variety of forms, it often faces challenges when applied to more complex systems, such as those involving ODEs/PDEs with discontinuous inputs or high-order derivatives. For example, current approaches struggle to accurately approximate discontinuities when identifying ODEs with a step input, which results in failures to capture the discontinuous terms (see Table \ref{tab:delta_func} in Section \ref{subsec:exp_discontinuous}). Another notable example is the identification of the Kuramoto–Sivashinsky equation $(u_t + u u_x + u_{xx} + u_{xxxx} = 0)$, where challenges arise from approximating the fourth-order derivative term $u_{xxxx}$. This challenge is further complicated by the presence of noise in measurements collected from real-world sensors (see Table \ref{tab:ks_result} in Section \ref{subsec:pde}), which can significantly undermine the accuracy of model discovery.

\subsection{Laplace-enhanced library}

Inspired by the use of the Laplace transform in various machine-learning tasks \cite{holt2022neural, cao2023lno}, we incorporate the Laplace transform into the construction of our library. This approach enables the transformation of candidate functions from the time domain into the Laplace domain, which provides a novel perspective for modeling dynamical systems, particularly those with complex features such as high-order derivatives or discontinuities. In this subsection, we denote the vector at time $t$ as $\mathtt{X}(t) = \mathtt{X}(d,k,n,t)\in\mathbb R^{\mathbf{d}}$, which can be represented as:
\begin{equation*}
\begin{split}
\mathtt{X}(t)=\left[
\mathtt{X}_1(t)\ \ \mathtt{X}_2(t)\ \ \cdots\ \ \mathtt{X}_{\mathbf{d}}(t)
\right]^\intercal
\end{split},
\end{equation*}
where $\mathtt{X}_i(t),\ (i=1, 2, \cdots, \mathbf{d})$ represents one element of $\mathtt{X}(t)$. 

To transform the candidate functions from the time domain to the Laplace domain, we apply the Laplace transformation to each element of $\mathtt{X}(t)$:
\begin{equation}\label{eq:laplace_rvs}
\begin{split}
    \mathbf{X}(s) & = \mathcal{L}\{\mathtt{X}(t)\}(s)\\
    & = \int_0^\infty e^{-st}\mathtt{X}(t)\mathrm{~d}t\\
    & = \left[\int_0^\infty e^{-st}\mathtt{X}_1(t)\mathrm{~d}t\ \ \int_0^\infty e^{-st}\mathtt{X}_2(t)\mathrm{~d}t\ \cdots\ \int_0^\infty e^{-st}\mathtt{X}_{\mathbf{d}}(t)\mathrm{~d}t\right]^\intercal\\
    & = \left[\mathbf{X}_1(s)\ \ \mathbf{X}_2(s)\ \ \cdots\ \ \mathbf{X}_{\mathbf{d}}(s)\right]^\intercal,
\end{split}
\end{equation}
where $\mathbf{X}(s) \in \mathbb{R}^{\mathbf{d}}$ denotes the vector of candidate functions in the Laplace domain, and each $\mathbf{X}_i(s)$ ($i = 1, 2, \cdots, \mathbf{d}$) represents the $i$-th element of the vector $\mathbf{X}(s)$.

In an ideal setting, the candidate functions enhanced by the Laplace transformation would integrate the variables from $t = 0$ to infinity. However, measurements collected from real-world sensors are always discrete and collected over a limited time interval, from $t_1$ to $t_m$. To approximate $\mathbf{X}(s)$ with a finite number of measurements, we employ the trapezoidal rule to compute a weighted sum of $\mathtt{X}(t_j)$ for $j = 1, 2, \dots, m$:
\begin{equation}\label{eq:laplace_trapezoidal}
    \mathbf{X}_i(s) = \int_0^\infty e^{-st}\mathtt{X}_i(t)\mathrm{~d}t \approx \sum_{j=1}^m e^{-s(t_j-t_1)}\mathtt{X}_i(t_j)\Delta t_j,
\end{equation}
where $\Delta t_j=t_{j+1} - t_j$ for $j=1,2,\cdots,m-1$, and $\Delta t_m = t_m - t_{m-1}$. For elements involving pure first-order to higher-order derivatives, we apply integration by parts within the Laplace transformation, followed by discretization using the trapezoidal rule:
\begin{equation}\label{eq:laplace_ibp}
\begin{split}
    % \mathcal{L}\left\{\frac{\mathrm{d}^k \boldsymbol u(t)}{\mathrm{d} t^k}\right\}(s) & = s^k \int_0^\infty e^{-st}\boldsymbol u(t)\mathrm{~d}t - \sum_{n=0}^{k-1}s^{k-n-1}\left.\frac{\mathrm{d}^n \boldsymbol u(t)}{\mathrm{d} t^n}\right|_{t=0}\\
    \mathcal{L}\left\{\frac{\mathrm{d}^k \boldsymbol u(t)}{\mathrm{d} t^k}\right\}(s) & = s^k \mathcal{L}\{\boldsymbol u(t)\}(s) - \sum_{n=0}^{k-1}s^{k-n-1}\left.\frac{\mathrm{d}^n \boldsymbol u(t)}{\mathrm{d} t^n}\right|_{t=0}\\
    & \approx s^k \sum_{j=1}^m e^{-s(t_j-t_1)}\boldsymbol u(t_j)\Delta t_j - \sum_{n=0}^{k-1}s^{k-n-1}\left.\frac{\mathrm{d}^n \boldsymbol u(t)}{\mathrm{d} t^n}\right|_{t=t_1}
\end{split}
\end{equation}
for $k=1,2,\dots, \mathbf k$. Here $\mathbf k$ represents the highest order derivative in the candidate function.

By carefully selecting different values of $s$ ($\{s \in \mathbb{C} \mid \text{Re}(s) > 0\}$) from $s_1$ to $s_L$, we derive various representations of $\mathbf{X}(s_l)$ for $l = 1, 2, \dots, L$. These representations contribute to the construction of the Laplace-enhanced library $\Theta:= \mathbb R^{\mathbf{d}}\times \mathbb C^{L}\rightarrow \mathbb{R}^{L \times \mathbf{d}}$

\begin{equation}\label{eq:build_lib}\centering
\begin{split}
\Theta(\bX) = \Theta(\bX, \bs) & =\left[
\mathbf{X}(s_1)\ \ \mathbf{X}(s_2)\ \ \cdots\ \ \mathbf{X}(s_L)\right]^\intercal\\
& =\left[\begin{array}{cccccccccc}
\mathbf{X}_1(s_1)   & \mathbf{X}_2(s_1)  && \cdots && \mathbf{X}_\mathbf{d}(s_1) \\
\mathbf{X}_1(s_2)   & \mathbf{X}_2(s_2)  && \cdots && \mathbf{X}_\mathbf{d}(s_2) \\
\vdots            & \vdots           && \ddots && \vdots  \\
\mathbf{X}_1(s_L)   & \mathbf{X}_2(s_L)  && \cdots && \mathbf{X}_{\mathbf{d}}(s_L) \\
\end{array}\right],\\
\end{split}
\end{equation}
where $\bs=\left[s_1, s_2, \cdots, s_L\right]^\intercal\in \mathbb C^L$, and each column within $\Theta(\bX)$ serves as a candidate function for the dynamical system. While the Laplace transformation typically involves complex frequencies, for simplicity, we commonly work with real numbers $s_l \in \mathbb{R}_{\geq 0}$. The choice of $s$ is critical and could impact performance, which we discuss its effects in Section \ref{subsec:exp_tri_and_hyper}. 

\subsection{Sparse regression}

As discussed earlier, the process of identifying a specific dynamical system involves constructing a Laplace-enhanced library with $\mathbf{d} = \frac{(kd+d+n+2)!}{n!(kd+d+2)!}$ candidate functions, based on $L$ measurements collected in the Laplace domain. Given that only a subset of candidate functions governs the complex dynamical system, we employ sparse regression \cite{brunton2016discovering} to extract the sparse coefficient vector $\xi\in\mathbb R^{\mathbf{d}}$, which identifies the active candidate functions in the library:
\begin{equation}\label{eq:sparse_regress}
	\Theta(\bX)\xi=\bm 0.
\end{equation}

It should be noted that a direct application of sparse regression in \eqref{eq:sparse_regress} may cause the weights to collapse to zero. To prevent this, a normal approach is to fix one weight to 1 at a time and conduct sparse regression iteratively across different fixed weights. We further denote $\xi_{-i}\in \mathbb R^{\mathbf d-1},\ (i=1, 2, \cdots, \mathbf{d})$ as a coefficient vector of $\xi$ except for the $i$-th element, and $\xi_{i}\in \mathbb R$ as the $i$-th element in $\xi$. The coefficient matrix $\Xi \in \mathbb{R}^{\mathbf d \times \mathbf d}$ comprises $\mathbf d$ vectors, each $\mathbf d$-dimensional. Each vector $\xi$ in $\Xi$ corresponds to an identified model and represents a governing equation. The diagonal elements of $\Xi$ are set to 1.0.
% \textcolor{red}{$\Xi = \left[\xi_{-1}, \xi_{-2}, \dots, \xi_{-\mathbf{d}}\right]^\intercal \in \mathbb{R}^{(\mathbf{d}-1)\times \mathbf{d}}$} as a coefficient matrix, 
$\mathbf{X}_{-i}(s)\in\mathbb R^{\mathbf{d}-1}$ denotes the vector $\mathbf{X}(s)$ without its $i$-th element, $\mathbf{X}_i(s)\in \mathbb R$. Specifically, we solve sparse regression problem for $i$ from $1$ to $\mathbf{d}$ iteratively:
\begin{equation}\label{eq:sparse_regress2}
	\Theta(\bX_{-i})\xi_{-i}=\Theta(\bX_{i}),
\end{equation}
where $\Theta(\bX_{-i})\in\mathbb R^{L\times(\mathbf{d}-1)}$ denotes the library $\Theta(\bX)$ except the $i$-th column, and $\Theta(\bX_{i})\in\mathbb R^{L}$ is the $i$-th candidate functions in $\Theta(\bX)$. 

To identify dynamical systems with sparse regression, the proposed LES-SINDy framework can employ various optimizers with their unique advantages. Stepwise Sparse Regression (SSR) \cite{mallows2000some} is a straightforward, greedy algorithm suitable for small datasets but may get stuck in local optima. Sparse Relaxed Regularized Regression (SR3) \cite{zheng2018unified} introduces flexibility through a relaxation parameter, which offers better convergence at the cost of increased complexity. Mixed-Integer Optimized Sparse Regression (MIOSR) \cite{bertsimas2016best} provides globally optimal solutions but is computationally expensive. Forward Regression Orthogonal Least Squares (FROLS) \cite{billings1996determination} reduce multicollinearity, which makes it ideal for more stable models. While Sequentially Thresholded Least Squares (STLS) \cite{daubechies2004iterative} offers a simple, iterative approach, it relies heavily on an appropriate selection of the threshold. Unless explicitly stated otherwise, STLS is employed as the default optimization method in our experiments. However, alternative optimizers are applied in our experiments where their specific advantages render them more suitable.

Once optimization is complete, the approach yields $\mathbf{d}$ distinct dynamical systems. We will then explore a suitable metric to identify the optimal dynamical system from among these candidates.

\subsection{Model evaluation}

To evaluate the accuracy of the dynamical systems identified by LES-SINDy, we first consider a metric inspired by log RMSE. This metric quantifies the discrepancy between the predicted trajectory of the identified model, $\tilde{\boldsymbol{u}}(t)\in \mathbb R^d$, and the actual measurements, $\boldsymbol{u}(t)$, over a specified time interval from $t = t_1$ to $t = t_m$. The metric $\epsilon_1:=\mathbb R^{\mathbf{d}}\times \mathbb R^{m}\times \mathbb R^{d\times m}\rightarrow\mathbb R$ can be represented as:
\begin{equation}\label{eq:log_rmse}
    \epsilon_1\left(\xi, \boldsymbol{t}, \boldsymbol{u}(\boldsymbol{t})\right) = \log \left[ {\left(\frac{1}{m} \sum_{j=1}^{m} \left\| \boldsymbol{u}(t_j) - \tilde{\boldsymbol{u}}(t_j, \xi) \right\|_2^2\right)}^{1/2} \right],
\end{equation}
where $\boldsymbol{t}=\left[t_1, t_2, \dots, t_m\right]\in \mathbb R^L_{\geq 0}$ denotes time vector, $\boldsymbol{u}(\boldsymbol{t})$ represents the actual measurement matrices from $t=t_1$ to $t_m$, and $\tilde{\boldsymbol{u}}(t_j; \xi)\in \mathbb R^{d}$ represents the predicted $d$ dimensional state from the identified model parameters $\xi$ at time $t_j$ ($j=1,2,\dots, m$). The state $\boldsymbol{u}(\boldsymbol{t})$ predicted by identified dynamical systems is typically solved using numerical methods, with Runge-Kutta methods \cite{butcher2016numerical} for ODEs and finite difference \cite{smith1985numerical, leveque2007finite} or spectral methods \cite{hesthaven2007spectral} for PDEs. The log RMSE metric provides a robust measure of model performance by emphasizing larger errors due to the logarithmic transformation while accounting for the different scales of the measurements. By inputting the identified model parameters $\xi$ and measurements over time, this metric outputs a scalar value that serves as an evaluation of the model's ability to accurately replicate the observed dynamics.

Another important metric we use to identify dynamical systems is AIC. While log RMSE provides a robust framework to evaluate identified dynamical systems, the process of identifying the most accurate and parsimonious model for nonlinear dynamical systems is challenging due to the need to balance model complexity and goodness of fit. AIC provides a statistical approach to this challenge by penalizing the number of free parameters while maximizing the likelihood of the model based on the given measurements. The AIC for each identified model, referred to as $\text{AIC}(\cdot, \cdot, \cdot)$, is computed as follows:
\begin{equation*}
    \text{AIC}\left(\xi, \boldsymbol{t}, \boldsymbol{u}(\boldsymbol{t})\right) = 2p + m \ln(2 \pi {\sigma}^2) + \frac{1}{{\sigma}^2} \sum_{j=1}^{m} (\boldsymbol{u}(t_i) - \tilde{\boldsymbol{u}}(t_j, \xi))^2,
\end{equation*}
where $p$ is the number of active elements in ${\xi}$, and ${\sigma}^2$ represents the estimated variance of the residuals. This metric evaluates the model’s fit to the measurements while imposing a penalty on models with more parameters to prevent overfitting. To further account for a limited number of measurements, the AIC can be modified using AICc, which is more accurate when $m$ is not large relative to $p$. The AICc is defined with $\epsilon_2:=\mathbb R^{\mathbf{d}}\times \mathbb R^{m}\times \mathbb R^{d\times m}\rightarrow\mathbb R$:
\begin{equation}\label{eq:aicc}
    \epsilon_2\left(\xi, \boldsymbol{t}, \boldsymbol{u}(\boldsymbol{t})\right) = \text{AIC}\left(\xi, \boldsymbol{t}, \boldsymbol{u}(\boldsymbol{t})\right) + \frac{2(p + 1)(p + 2)}{m - p - 2}.
\end{equation}
By using AICc, we systematically reduce a large set of candidate functions, which ensures the output model is both accurate and parsimonious. This approach effectively reduces the risk of overfitting and increases the robustness of the identified models, which makes LES-SINDy a more dependable tool for discovering governing equations in nonlinear dynamical systems. This strategy effectively mitigates the risk of overfitting and improves the robustness of identified models, which makes LES-SINDy a more reliable tool for discovering governing equations in nonlinear dynamical systems.

% \textcolor{red}{AICc is a relative score, which cannot compare across different tasks}

\begin{algorithm}
\caption{Laplace-Enhanced Sparse Identification of Nonlinear Dynamical Systems. Given time-series measurements, the algorithm outputs identified models to characterize the underlying dynamical systems.}\label{alg:les-sindy}
\textbf{Input} time stamp: $\boldsymbol{t}=\left[t_1, t_2, \cdots, t_m\right]^\intercal$, $t_1,t_2,\cdots,t_m \in\mathbb R_{\geq0}$;\\
\textbf{Input} complex frequency: $\bs=\left[s_1, s_2, \cdots, s_L\right]^\intercal\in \mathbb C^L$, $\text{Re}(s_l) > 0$ for $l=1,2,\dots, L$;\\
\textbf{Input} measurements: $\boldsymbol u(t)=\left[u_1(t),u_2(t),\cdots,u_d(t)\right]^\intercal\in \mathbb R^{d}$, $t=t_1,t_2,\cdots,t_m$;\\
\textbf{Input} order of (time) derivatives: $k$; number of tensor product copies: $n$. \vspace{0.05 in}
% {\textbf{Require}} Approximate samples from the posterior distribution $p(\beta \mid D)$
\begin{algorithmic}[1]
\State \textcolor{darkred}{\textbf{Build Laplace-Enhanced Library}\Comment{Step 1}}
\For{$i = 1$ to $m$}
    \State Construct candidate function in the time domain: $\mathtt{X}(d,k,n,t_i)=\bigotimes^n \mathcal{D}^{k}\bu(t_i)$
\EndFor
\For{$i = 1$ to $L$}
    \State Construct candidate function $\mathbf{X}(s_i)$ in the Laplace domain following \eqref{eq:laplace_rvs}-\eqref{eq:laplace_ibp}
\EndFor
\State Build a library $\Theta(\bX, \bs)$ following \eqref{eq:build_lib} \vspace{0.05 in}

% \noindent{\NoNumber \textbf{Sparse Regression}}
\State \textcolor{darkgreen}{\textbf{Sparse Regression}\Comment{Step 2}}
\For{$i = 1$ to $\mathbf d$}
    \State Set $\xi_{i}=1.0$ 
    \State Find $\xi_{-i}\in \mathbb R^{\mathbf d-1}$ by solving \eqref{eq:sparse_regress2}
    \State Save the identified model: $\xi^{(i)} = \xi$
\EndFor \vspace{0.05 in}
\State \textcolor{darkblue}{\textbf{Model Evaluation}\Comment{Step 3}}
\State Construct a vector to store AICc: $\boldsymbol\epsilon=\left[\epsilon^{(1)}, \epsilon^{(2)},\cdots,\epsilon^{(\mathbf d)}\right]^\intercal$
\For{$i = 1$ to $\mathbf d$}
    \State Calculate $\epsilon^{(i)} = \epsilon_2\left(\xi^{(i)}, \boldsymbol{t}, \boldsymbol{u}(\boldsymbol{t})\right)$ following \eqref{eq:aicc}
\EndFor \vspace{0.05 in}
\end{algorithmic}
{\textbf{Output}} Identified model with the lowest AICc: $\argmin_{\xi}\boldsymbol\epsilon$
\end{algorithm}

Comparing the metrics as mentioned earlier, the log RMSE measures prediction accuracy and is ideal for comparing models with similar complexity. We will focus on using this metric to optimize hyperparameter selection across our experiments. On the other hand, AICc, which balances the evaluation of model fit and complexity, helps to prevent overfitting. In our experiments, we use AICc to compare identified models and to benchmark our results against previous methods like SINDy \cite{brunton2016discovering, rudy2017data} and Weak-SINDy \cite{messenger2021weak} when hyperparameters are held constant.

To sum up, LES-SINDy is an algorithm designed to identify governing equations for complex dynamical systems using sparse regression techniques. The algorithm begins with the collection of time-series measurements from the system under study. These measurements are then transformed into the Laplace domain across various complex frequencies, which enables the construction of a comprehensive library of candidate functions. The algorithm applies sparse regression techniques to update multiple candidate models, which are further evaluated through AICc to form a parsimonious model. Figure \ref{fig:framework} and Algorithm \ref{alg:les-sindy} illustrate the step-by-step process of the LES-SINDy framework, which demonstrates how time-series measurements are applied to capture the essential behavior of the underlying dynamical system.

% \clearpage
\section*{Data Availability}
Datasets used in this paper will be publicly available after publication. Data can be downloaded from \href{http://github.com/haoyangzheng1996/LES-SINDy/tree/main/data}{github.com/haoyangzheng1996/les-sindy/tree/main/data}.

\section*{Code Availability}

Code will be released at \href{http://github.com/haoyangzheng1996/LES-SINDy}{github.com/haoyangzheng1996/les-sindy} after publication.

% \clearpage

\bibliography{refs}

\clearpage

\begin{appendices}

\section{Candidate function construction}\label{appendix:tensor_candidate}

% \textcolor{red}{\textbf{Require Modification.}}

As mentioned in Section \ref{subsec:tensor_product}, we define the $n$-th tensor product of $\mathcal D$ to capture higher-order interactions:

$$
\bigotimes^n \mathcal D = \underbrace{\mathcal D \otimes \mathcal D \otimes \cdots \otimes \mathcal D}_{n \text{ times}}.
$$

This operation generates all possible products of $n$ elements from $\mathcal D$, which includes repeated elements. The resulting set includes monomials of degree $n$ formed by the original functions and their derivatives. The output dimension $\mathbf{d}$ of the extended candidate function set is given by the combinatorial formula  $\mathbf{d}=\frac{({kd+d+n+1})!}{n!({kd+d+1})!}$. This formula calculates the number of ways to choose $n$ elements from $kd+d+2$ options with replacement.

To further clarify how the extended set is constructed, consider the following examples:

\begin{enumerate}
    \item When $n = 1$ and $k = 0$, the basis functions are simply:
       $$
       \left\{1, \ t, \ \boldsymbol{u}\right\},
       $$
       where $\boldsymbol{u} = [u_1(t), \ u_2(t), \ \dots, \ u_d(t)]^\intercal$. This set includes the constant term, time $t$, and the system variables $u_i(t)$.
    \item When $n = 1$ and $k = 1$, the basis functions expand to include the first derivatives:
       $$
       \left\{1, \ t, \ \boldsymbol{u}, \ \frac{\mathrm{d} \boldsymbol{u}}{\mathrm{d} t}\right\}.
       $$
       This accounts for both the variables and their rates of change and captures linear dynamics.
    \item When $n = 1$ and $k = 2$, the basis functions further include second derivatives:
       $$
       \left\{1, \ t, \ \boldsymbol{u}, \ \frac{\mathrm{d} \boldsymbol{u}}{\mathrm{d} t}, \ \frac{\mathrm{d}^2 \boldsymbol{u}}{\mathrm{d} t^2}\right\}.
       $$
       This set is capable of representing acceleration and other second-order effects in the system.
   \item When $n = 2$ and $k = 1$, the basis functions include all products of two elements from $\mathcal D$:
       $$
       \left\{1, \ t, \ \boldsymbol{u}, \ \frac{\mathrm{d} \boldsymbol{u}}{\mathrm{d} t}, \ t^2, \ t \boldsymbol{u}, \ t \frac{\mathrm{d} \boldsymbol{u}}{\mathrm{d} t}, \ \boldsymbol{u}^2, \ \boldsymbol{u} \frac{\mathrm{d} \boldsymbol{u}}{\mathrm{d} t}, \ \left(\frac{\mathrm{d} \boldsymbol{u}}{\mathrm{d} t}\right)^2\right\}.
       $$
       This set captures quadratic terms and interactions between variables and their derivatives.
    \item When $n = 2$ and $k = 2$, a total number of 15 basis functions are given as follows:
    $$
       \left\{1, \ t, \ \boldsymbol{u}, \ \frac{\mathrm{d} \boldsymbol{u}}{\mathrm{d} t}, \ \frac{\mathrm{d}^2 \boldsymbol{u}}{\mathrm{d} t^2}, \ t^2, \ t \boldsymbol{u}, \ t \frac{\mathrm{d} \boldsymbol{u}}{\mathrm{d} t}, 
       \ t \frac{\mathrm{d}^2 \boldsymbol{u}}{\mathrm{d} t^2}, \ \boldsymbol{u}^2, \ \boldsymbol{u} \frac{\mathrm{d} \boldsymbol{u}}{\mathrm{d} t}, \ \boldsymbol{u} \frac{\mathrm{d}^2 \boldsymbol{u}}{\mathrm{d} t^2}, \ \left(\frac{\mathrm{d} \boldsymbol{u}}{\mathrm{d} t}\right)^2, \ \frac{\mathrm{d} \boldsymbol{u}}{\mathrm{d} t}\frac{\mathrm{d}^2 \boldsymbol{u}}{\mathrm{d} t^2},\ \left(\frac{\mathrm{d}^2 \boldsymbol{u}}{\mathrm{d} t^2}\right)^2\right\}.
       $$
\end{enumerate}

\clearpage

\section{High-order ODEs}\label{appendix:high-order}
We consider a fourth-order linear ODE given by:
$$u_{tttt}+8u_{tt}+16u=0,$$
with the initial conditions $u(0)=u_{t}(0)=u_{tt}(0)=0, u_{ttt}(0)=1.$ The time interval for the simulation runs from 0 to 20, with 200 uniformly distributed time points. The system is free of noise, which allows us to focus purely on the model identification performance. For model identification, we apply the thresholding least squares method with a threshold value of 0.01. 

The results of this experiment are compared with the standard SINDy:
\begin{table}[htbp]
\centering
\begin{tabular}{cccccccc}\toprule
Results    & 1      & 2     & 3     & 4     & 5     & 6       & 7      \\ \midrule
$u_{tttt}$    & 1      & 0.516 & 0.125 & 0.129 & 0.062 & -2.324  & 0.263  \\
$u_{ttt}$       & 0.066  & 1     &       & 0.248 &       & 0.301   &        \\
$u_{tt}$        & 7.950   & 3.963 & 1     & 0.954 & 0.496   & -16.334 & 1.759  \\
$u_{t}$         & 0.264  & 3.996 &       & 1     &       & 1.128   & 0.039  \\
$u $         & 15.843 & 8.207 & 1.991 & 1.901 & 1     & -27.711 & 2.823  \\
$t$          &        &       &       &       &       & 1       & -0.079 \\
$1$          &        &       &       &       &       & -9.510   & 1      \\ \midrule
AICc         & -587.8 & 268.5 & -845.9 & 269.7 & -595.5 & 597.9 & 383.6 \\ \bottomrule
\end{tabular}
\end{table}

Results from the proposed LES-SINDy:
\begin{table}[htbp]\centering
\begin{tabular}{cccccccc}\toprule
Results     & \textbf{1}         & 2        & 3          & 4        & 5     & 6        & 7        \\ \midrule
$u_{tttt}$       & \textbf{1}         & 1369.330  & 0.125      & -67.779  & 0.062 & 1786.555 & -510.172 \\
$u_{ttt}$          &           & 1        &            & -0.022   &       & 0.702    & -0.204   \\
$u_{tt}$           & \textbf{8.000}     & 10970.850 & 1          & -542.913 & 0.500   & 14310.84 & -4086.530 \\
$u_{t}$            &           & -23.204  &            & 1        &       & -27.015  & 7.683    \\
$u$             & \textbf{16.000}    & 21975.610 & 2.000      & -1087.050 & 1     & 28656    & -8182.750 \\
$t$             &           & 0.984    &            & -0.036   &       & 1        & -0.283   \\
$1$             &           & -3.402   &            & 0.127    &       & -3.535   & 1       \\ \midrule
AICc         &  \textbf{-2559.3}& -1327.9 & -1722.8 & -1421.7 & -996.7 & -1550.5 & -1400.7 \\ \bottomrule
\end{tabular}
\end{table}

\clearpage

\section{ODEs with discontinuous inputs}

\subsection{ODEs with a step function}

This experiment focuses on a first-order ODE with a step input at $t_0$, given by:
$$a u_{t} + b u + cH(t-t_0)=0,$$
where $a=1.0$, $b=2.0$, and $c=1.0$. The step function $H(t-t_0)$ introduces a discontinuity at $t_0$, which makes this system an ideal test case for evaluating the performance of model identification methods in the presence of abrupt changes.

We first examine the results of both SINDy and LES-SINDy under noise-free conditions to assess their respective abilities to capture the system dynamics. The results for SINDy are given as follows:
\begin{table}[htbp]
\centering
\begin{tabular}{cccccc}
\toprule
Results       &       1      & 2    & 3    & {4}      \\ \midrule
$\delta(t)$            & 1       & -0.505 & 0.252  & 0.505             \\
$H(t)$                 & -1.980  & 1      & -0.500 & -1.000            \\
$u$                    & 3.960   & -2.000 & 1      & 2.000             \\
$u_{t}$                   & 1.980   & -1.000 & 0.500  & 1                 \\ \midrule
AICc       & 423.2 & 423.2 & 423.2  & 464.1 \\ \bottomrule
\end{tabular}
\end{table}

Here are the results yielded by LES-SINDy:
\begin{table}[htbp]
\centering
\begin{tabular}{cccccc}
\toprule
Results       &       1      & 2    & \textbf{3}    & {4}      \\ \midrule
$\delta(t)$   &  1    &      &    &    \\
$H(t)$        &  691.188     & {1}       &  \textbf{0.500}   & 1.001   \\
$u$           & -1385.021    & {2.002}   & \textbf{1}    & 2.005   \\
$u_{t}$          & -685.952     & {0.999}   & \textbf{0.500}   & 1    \\ \midrule
AICc       & -1198.1 & -1468.1 & \textbf{-1468.2}  & -1468.0 \\ \bottomrule
\end{tabular}
\end{table}

Additionally, we further analyze the performance of LES-SINDy as the noise level is progressively increased from 0\% to 10\%, which provides insight into the robustness of the proposed method under noisy conditions.
\begin{table}[!htbp]
\centering\caption{Results of LES-SINDy applied to ODE with a step input with varying levels of noise.}\label{appendix:tab:step_noise}
\begin{tabular}{cc|ccccccc}
\toprule
\multicolumn{2}{c|}{\multirow{2}{*}{results}}  &       \multicolumn{6}{c}{candidate functions} \\ 
\multicolumn{2}{c|}{}                          & {$u_t$}  & {$u$}     & $t$      & $1$     & $\delta (t-t_0)$  & {$H(t-t_0)$}  \\ \midrule
\multirow{12}{*}{\rotatebox[origin=c]{90}{noise level}} 
              & truth & 1.0   & 2.0   &        &       &        & 1.0   \\ \midrule
              & clean & 1.000 & 2.001 &        &       &        & 1.000 \\
              & 1\%   & 1.000 & 2.002 &        &       &        & 1.001 \\
              & 2\%   & 1.000 & 1.895 &        &       &        & 0.968 \\
              & 3\%   & 1.000 & 1.646 &        &       &        & 0.890 \\
              & 4\%   & 1.000 & 1.251 &        &       &        & 0.768 \\
              & 5\%   & 1.000 & 1.232 &        &       &        & 0.761 \\
              & 6\%   & 1.000 & 1.426 & -0.018 &       & -0.034 & 0.816 \\
              & 7\%   & 1.000 & 1.475 & -0.025 &       & -0.048 & 0.830 \\
              & 8\%   & 1.000 & 1.082 & -0.013 &       & -0.012 & 0.696 \\
              & 9\%   & 1.000 & 0.295 & -0.004 &       & 0.048  & 0.390  \\ 
              & 10\%  & 1.000 & 0.730 & -0.019 &       & 0.023  & 0.518  \\ \bottomrule
\end{tabular}
\end{table}

\clearpage

\subsection{ODEs with a delta function}

We consider a second-order ODE with a delta input at $t_0$, which is described by:
$$au_{tt} + bu_{t} + c u + d\delta(t-t_0) = 0$$
with $a=1.0$, $b=4.0$, $c=4.0$, and $d=1.0$. The Dirac delta function $\delta(t-t_0)$ represents an instantaneous impulse at time $t_0$, which makes this system ideal for testing the performance of identification algorithms when handling abrupt, short-term inputs.

We first present results from both SINDy and LES-SINDy under noise-free conditions to evaluate how well each method captures the system dynamics. Here are the results of SINDy:
\begin{table}[htbp]
\centering
\begin{tabular}{cccccc}
\toprule
Results         &    1      & 2          &    3    & 4       & {5}                \\ \midrule
$\delta(t)$            & 1                   & 53.000              & -24.999 & -24.998 & -99.980    \\
$H(t)$                 &                     & 1                   &         &         &     \\
$u$                    & -0.040              & -10.141             & 1       & 1.000   & 3.999     \\
$u_{t}$                   & -0.040              & -2.109              & 1.000   & 1       & 3.999      \\
$u_{tt}$                  & -0.010              & -0.539              & 0.250   & 0.250   & 1.000      \\ \midrule
AICc            &  126.8 & 82.1 & 126.8 & 126.8 & 126.9 \\ \bottomrule\end{tabular}
\end{table}

Here are the results from LES-SINDy:
\begin{table}[htbp]
\centering
\begin{tabular}{cccccc}
\toprule
Results         &    1      & 2          &    \textbf{3}    & 4       & 5                \\ \midrule
$\delta(t)$     & {1}       & -61.457    & \textbf{0.250}   &  0.250  & 1.000   \\
$H(t)$          &           & 1          &         &         &                  \\
$u$             & {4.000}   &  -222.920  & \textbf{1}       & 1.000   &  3.999  \\
$u_{t}$            & {4.000}   &  -269.314  & \textbf{1.000}   & 1       & 3.999   \\
$u_{tt}$           & {1.000}   &  -61.259   & \textbf{0.250}   &  0.250  & 1       \\ \midrule
AICc            &  -1912.7 & -195.3 & \textbf{-1921.4} & -1914.4 & -1913.3 \\ \bottomrule
\end{tabular}
\end{table}

Additionally, we extend the analysis by applying LES-SINDy to datasets with increasing levels of measurement noise, ranging from 0\% to 20\%, to assess its robustness and accuracy in the presence of noisy inputs.
\begin{table}[!htbp]
\centering\caption{Results of LES-SINDy applied to ODE with a delta input with varying levels of noise.}\label{appendix:tab:delta_noise}
\begin{tabular}{cc|ccccccccc}
\toprule
\multicolumn{2}{c|}{\multirow{2}{*}{results}}    &      \multicolumn{7}{c}{candidate functions} \\ 
\multicolumn{2}{c|}{}                            & {$u_{tt}$} & {$u_t$}  & {$u$}     & $t$      & $1$     & {$\delta (t-t_0)$}  & {$H(t-t_0)$}  \\ \midrule
\multirow{22}{*}{\rotatebox[origin=c]{90}{noise level}} 
   & truth & 1.0   & 4.0   &  4.0      &       &        & 1.0   \\ \midrule
   & clean & 1.000  & 4.000 & 4.000 &       &   & 1.000 &        \\
   & 1\%   & 1.000  & 4.000 & 3.999 &       &   & 1.001 &        \\
   & 2\%   & 1.000  & 3.996 & 4.011 &       &   & 0.989 &        \\
   & 3\%   & 1.000  & 4.001 & 3.997 &       &   & 1.001 &        \\
   & 4\%   & 1.000  & 3.997 & 4.008 &       &   & 0.990 &        \\
   & 5\%   & 1.000  & 4.032 & 3.902 &       &   & 1.103 &        \\
   & 6\%   & 1.000  & 4.005 & 3.984 &       &   & 1.017 &        \\
   & 7\%   & 1.000  & 3.995 & 4.001 &       &   & 1.000 &        \\
   & 8\%   & 1.000  & 4.006 & 3.972 &       &   & 1.037 &        \\
   & 9\%   & 1.000  & 4.004 & 3.968 &       &   & 1.077 &        \\
   & 10\%  & 1.000  & 3.916 & 4.232 &       &   & 0.851 &        \\
   & 11\%  & 1.000  & 3.866 & 4.369 &       &   & 0.787 &        \\
   & 12\%  & 1.000  & 3.943 & 4.252 &       &   & 0.836 & -0.003 \\
   & 13\%  & 1.000  & 3.898 & 4.294 &       &   & 0.792 &        \\
   & 14\%  & 1.000  & 3.896 & 4.270 &       &   & 0.897 &        \\
   & 15\%  & 1.000  & 3.656 & 4.840 & 0.014 &   & 0.143 & -0.061 \\
   & 16\%  & 1.000  & 3.865 & 4.367 &       &   & 0.809 &        \\
   & 17\%  & 1.000  & 3.755 & 4.650 & 0.013 &   & 0.294 & -0.061 \\
   & 18\%  & 1.000  & 3.893 & 4.227 & 0.005 &   & 0.826 & -0.022 \\
   & 19\%  & 1.000  & 3.808 & 4.552 & 0.006 &   & 0.577 & -0.028 \\
   & 20\%  & 1.000  & 3.793 & 4.585 & 0.008 &   & 0.448 & -0.039   \\ \bottomrule
\end{tabular}
\end{table}

\clearpage

\section{ODEs with trigonometric and hyperbolic functions}
In this series of experiments, we examine the ability of LES-SINDy to accurately identify ODEs that involve trigonometric and hyperbolic functions. Each experiment is designed to test LES-SINDy’s performance with different types of oscillatory and exponential behavior.

\subsection{Learn the sine function}
We consider the following second-order ODE with a sine input:
$$u_{tt} + 15u + 2sin(3t)=0,$$ 
with initial conditions $u(0)=0$ and $u_{t}(0)=0$. This experiment tests the algorithm’s ability to learn periodic dynamics driven by a sine function. We present the results from LES-SINDy as follows:

\begin{table}[htbp]
\centering
\begin{tabular}{cccccc}
\toprule
Results    & \textbf{1}      & 2      & \textbf{3}      & 4\\ \midrule
$u_{tt}$     & \textbf{1}             & -225.100  & 0.067      & -17240000  &          \\ 
$u_{t}$      &               & 1         &            & 62.250         &          \\
$u$       & \textbf{15.000}        & -3375.000 & 1          & -258700000 &          \\
$t$       &               &           &            & 1              &          \\
$1$       &               &           &            & -3.141         &          \\
$sin(t)$  & \textbf{2.000}         & -449.400  & 0.133      & -34490000  &          \\
$cos(t)$  &               & 0.127     &            & 20.310         &          \\
$sinh(t)$ &               & 0.533     &            & 45.700         &          \\
$cosh(t)$ &               & -0.533    &            & -45.700        &          \\ \midrule
AICc      & \textbf{-2456.6} & -1424.4 & -1139.4 & -1654.1 \\ \bottomrule
        &               &           &            &                &          \\ \\ \toprule
Results    & 5         & \textbf{6}      & 7      & 8 & 9\\ \midrule
$u_{tt}$     & 6674000   & \textbf{0.500}     & 14360  &                & 129.200  \\ 
$u_{t}$      & -3.620    &           & 4.469      & 1.872          & -40.240  \\
$u$       & 100100000 & \textbf{7.500}     & 215400 & 5.850          & 1834.000 \\
$t$       & -0.368    &           &            &                &          \\
$1$       & 1         &           &            &                &          \\
$sin(3t)$  & 13350000  & \textbf{1}         & 28720  & 1.720          & 238.700  \\
$cos(3t)$  & -4.371    &           & 1          & 0.235          & -0.235   \\
$sinh(3t)$ & -11.500   &           & 3.829      & 1              & -1.000   \\
$cosh(3t)$ & 11.500    &           & -3.829     & -1.000         & 1        \\ \midrule
AICc      & -1798.4 & \textbf{-2267.5} & -1944.7 & 69918.8 & -491.4 \\ \bottomrule
\end{tabular}
\end{table}

\clearpage

\subsection{Learn the cosine function}
Next, we investigate the identification of an ODE with a cosine input:
$$u_{tt} + 4u + cos(t)=0,$$ 
with initial conditions $u(0)=0$ and $u_{t}(0)=0$. This test evaluates the model's capability to capture oscillations from a cosine source. Results from LES-SINDy are provided here:

\begin{table}[htbp]
\centering
\begin{tabular}{cccccc}
\toprule
Results    & \textbf{1}      & 2      & \textbf{3}      & 4         \\ \midrule
$u_{tt}$     & \textbf{1}                & -3849000000      & 0.250 & 13100000     &            \\
$u_{t}$      &                  & 1                &       & 28.380       &            \\
$u$       & \textbf{4.000}            & -15400000000     & 1     & 52400000     &            \\
$t$       &                  & 11.500           &       & 1            &            \\
$1$       &                  & 54.850           &       & -4.252       &            \\
$sin(t)$  &                  & 434.900          &       & 2.728        &            \\
$cos(t)$  & \textbf{1.000}            & -3849000000      & 0.250 & 13100000     &            \\
$sinh(t)$ &                  & -3.822           &       & 23.750       &            \\
$cosh(t)$ &                  & 3.822            &       & -23.750      &            \\ \midrule
AICc      & \textbf{-1891.3} & -1406.4 & -1146.1 & -1512.8  \\ \bottomrule
        &               &           &            &                &          \\ \\ \toprule
Results     & 5         & {6}      & \textbf{7}      & 8 & 9\\ \midrule
$u_{tt}$     & -23070000000     & -7363000         & \textbf{1.000} & 1152.000     & -3040      \\
$u_{t}$      & -27130           & 6.582            &       & 41.990       & 187.300    \\
$u$       & -92290000000     & -29450000        & \textbf{4.000} & 4608.000     & -12030     \\
$t$       & 161.900          & 0.337            &       &              &            \\
$1$       & 1                & -0.999           &       &              & 0.121      \\
$sin(t)$  & 1556.000         & 1                &       & 0.884        &            \\
$cos(t)$  & -23070000000     & -7363000         & \textbf{1}     & 1152.000     & -2944.000  \\
$sinh(t)$ & -4.055           & 4.784            &       & 1            & -1.000     \\
$cosh(t)$ & 4.055            & -4.784           &       & -1.000       & 1          \\  \midrule
AICc      & -1580.9 & -1675.2 & \textbf{-1751.0} & -549.4 & -686.6 \\    \bottomrule
\end{tabular}
\end{table}

\clearpage

\subsection{Learn the hyperbolic sine function}
We now consider an ODE involving the hyperbolic sine function:
$$u_{tt} + 4u + sinh(2t)=0,$$
with initial conditions $ u(0)=0$ and $u_{t}(0)=0$. This case is designed to assess the model’s performance in identifying exponentially growing behavior. The results obtained from LES-SINDy are presented here:

\begin{table}[htbp]
\centering
\begin{tabular}{cccccccccc}\toprule
Results      & \textbf{1}      & 2      & 3      & 4       & 5      & 6       & 7        & \textbf{8}     & 9     \\ \midrule
$u_{tt}$      & \textbf{1}      & 1.443  & -2.586 & -4.046  & 0.304  & -15.500   & 3.319  & \textbf{1.000} & 0.750  \\
$u_{t}$       &                 & 1      & 0.708  & -3.151  & 0.648  & -5.792   & 1.996   &                & 0.726 \\
$u$        & \textbf{4.000}  & 5.377  & 1      & -12.953 & 2.535  & -45.232  & 12.515  & \textbf{4.000} & 4.500   \\
$t$        &                 & -0.084 & -0.203 & 1       & -0.134 & 1.524    & -0.649  &                & -0.220 \\
$1$        &                 & 1.558  & 1.371  & -4.335  & 1      & -11.571  & 3.836   &                & 1.319 \\
$\sin(2t)$  &                 & -0.123 & -0.122 & 0.325   & -0.085 & 1        & -0.333  &                & -0.110 \\
$\cos(2t)$  &                 & 0.355  & 0.353  & -0.915  & 0.251  & -2.984   & 1       &                & 0.326 \\
$\sinh(2t)$ & \textbf{1.000}  & 1.416  & -2.578 & -3.606  & 0.174  & -15.203  & 2.977   & \textbf{1}     & 1.000     \\
$\cosh(2t)$ &                 & 1.212  & 1.019  & -3.446  & 0.755  & -8.653   & 2.857   &                & 1     \\ \midrule
AICc      &    \textbf{-1050.5}   & 741.8 & 571.8 & 584.5 & 666.1 & 571.7 & 608.5 & \textbf{-910.5} & 599.5 \\ \bottomrule
\end{tabular}
\end{table}

\subsection{Learn the hyperbolic cosine function}
Finally, we examine an ODE driven by the hyperbolic cosine function:
$$u_{tt} - 4u + cosh(2t)=0,$$ 
with initial conditions $u(0)=0$ and $u_{t}(0)=0$. This experiment tests LES-SINDy’s accuracy in identifying dynamics involving hyperbolic growth. The results from LES-SINDy are detailed as follows:

\begin{table}[htbp]
\centering
\begin{tabular}{cccccccccc}\toprule
Results     & \textbf{1}        & 2         & {3}        & 4         & 5         & 6       & 7       & 8      & \textbf{9}  \\ \midrule
$u_{tt}$     & \textbf{1}        & 0.688     & {-0.250}   & 0.281     & -0.500    & -0.469  & 2.500     & 1.125  & \textbf{1.000}  \\
$u_{t}$      &                   & 1         &                   & 4.251     & -2.285    & -13.932 & -77.156 & 1.395  &    \\
$u$       & \textbf{-4.000}   & -4.500    & {1}        & -9.500    & 6.000     & 29.625  & 145.000     & -7.000     & \textbf{-4.000} \\
$t$       &                   & 0.254     &                   & 1         & -0.582    & -3.898  & -22.578 & 0.369  &    \\
$1$       &                   & -0.436    &                   & -1.718    & 1         & 6.359   & 36.438  & -0.634 &    \\
$sin(2t)$  &                   & -0.070    &                   & -0.278    & 0.162     & 1       & 5.686   & -0.102 &    \\
$cos(2t)$  &                   &           &                   &           &           & 0.135   & 1       &        &    \\
$sinh(2t)$ &                   & 0.688     &                   & 2.709     & -1.577    & -9.694  & -55.000     & 1      &    \\
$cosh(2t)$ & \textbf{1.000}    & 0.312     & {-0.250}   & -1.375    & 0.750     & 5.656   & 37.875  & 0.250   & \textbf{1  }\\ \midrule
AICc   & \textbf{-2369.4} & 71.6 & -1629.3 & 325.8 & 266.9 & 376.5 & 5025.1 & 104.9 & \textbf{-2089.3}       \\ \bottomrule
\end{tabular}
\end{table}

\clearpage

\section{Nonlinear ODE systems}

In this section, we evaluate LES-SINDy’s ability to identify two well-known nonlinear ODE systems: the Lorenz system and the Lotka-Volterra model. These experiments aim to test the method's performance in capturing the complex, nonlinear interactions characteristic of such systems.

\subsection{The Lorenz systems}

The Lorenz System exhibits chaotic behavior that is highly sensitive to initial conditions, which makes it an ideal candidate for testing model identification techniques. We present the results from LES-SINDy, which demonstrates its ability to capture the underlying dynamics of the Lorenz system.

\begin{table}[htbp]
\caption{Discovered dynamical systems within the Lorenz system using LES-SINDy. The table shows candidate basis functions in the first column. Results 3, 5, and 6, highlighted in bold, exhibit small error bars and accurately capture the true dynamics of the Lorenz system.}\label{tab:lorenz_appendix}\centering
\begin{tabular}{cccccccccc}\toprule
results & 1         & 2        & \textbf{3}          & 4          & \textbf{5}          & \textbf{6}           \\ \midrule
$x$   & 1         & -1.085   &                    & -59.768    & \textbf{-23.974}    & \textbf{-3.180}      \\
$y$   & -2.037    & 1        &                    & -7.517     & \textbf{-3.015}     & \textbf{3.179}       \\
$z$   & -0.320    &          & \textbf{1}          &            &            & \textbf{2.667}       \\
$\dot x$  &           &          &                    & 1          & \textbf{0.401}      & \textbf{-0.318}      \\
$\dot y$  &           &          &                    & 2.493      & \textbf{1}          &             \\
$\dot z$  & 0.240     &          & \textbf{0.375}      &            &            & \textbf{1}           \\
$tx$  &           &          &                    &            &            &             \\
$ty$  &           &          &                    &            &            &             \\
$tz$  &           &          &                    &            &            &             \\
$xy$  &           &          & \textbf{-0.375}     &            &            & \textbf{-1.000}      \\
$xz$  &           &          &                    & 2.492      & \textbf{1}          &             \\
$yz$  &           &          &                    &            &            &             \\ \bottomrule
        &               &           &            &                &          \\ \\ \toprule
    % &           &          &                    &            &            &             \\ \hline\hline
results & 7         & 8        & 9                  & 10         & 11         & 12          \\ \midrule
$x$   & -1411.677 & -354.255 & -18388.000 & 50935.418  & 74642.000  & -160216.000 \\
$y$   & 1411.661  & 264.425  & 19032.000  & -50923.818 & -74716.000 & 160080.000  \\
$z$   &           & -9.860   & -43.188    & -3.626     & -5.000     & 1.000           \\
$\dot x$  & -141.089  & -27.121  & -1893.500  & 5094.602   & 7520.750   & -15896.000  \\
$\dot y$  &           & 3.505    & -25.625    & -4.560     & -2.497     & 12.000      \\
$\dot z$  &           & -3.535   & 5.563      & -0.855     & -1.500     & -5.000      \\
$tx$  & 1         & -1.001   & -153.210   & 0.124      &            & 17.391      \\
$ty$  & -0.998    & 1        & 155.745    & -0.116     &            & -17.688     \\
$tz$  &           &          & 1          &            &            & -0.109      \\
$xy$  &           & 3.671    & 4.750      & 1          & 1.250      & 2.500       \\
$xz$  &           & 3.372    & -18.750    & -2.595     & 1          & 5.000       \\
$yz$  &           &          &            &            &            & 1           \\ \bottomrule
\end{tabular}
\end{table}

\clearpage

\subsection{The Lotka-Volterra model}
The Lotka-Volterra model, also known as the predator-prey model, describes the interaction between two species: one as prey and the other as predator. This experiment assesses LES-SINDy’s ability to capture nonlinear population dynamics. The identified dynamical systems from LES-SINDy are provided as follows to evaluate their performance in identifying this model.

\begin{table}[!htbp]
\caption{Discovered dynamical systems within the Lotka-Volterra model using LES-SINDy. The table shows candidate basis functions in the first column. Results 3, 5, and 6, highlighted in bold, exhibit small error bars and accurately capture the true dynamics of the Lotka-Volterra model.}\label{tab:lv_appendix}\centering
    \begin{tabular}{cccccccc}
        \toprule
        Results    & 1 & 2 & 3 & 4 & \textbf{5} \\ \midrule
                $1$  & 1 & -2.337 &   &   &    \\ 
        $t$  & -0.182 & 1 &   &   &     \\ 
        $x$  & -153.363 & 0.865 & 1 & 0.195 & \textbf{-1.004} \\ 
        $y$  & -531.305 & 1.781 & 5.052 & 1 & \textbf{-5.071}  \\ 
        $\dot x$  & 152.071 & -0.326 & -0.996 & -0.194 & \textbf{1} \\ 
        $\dot y$  & -530.992 & 1.555 & 5.052 & 1.000 & \textbf{-5.071} \\ 
        $tt$  &  & -0.120 &   &   &   &    \\ 
        $tx$  &   &   &   &   &  \\ 
        $ty$  & & 0.134 &   &   &   \\ 
        $xy$  & 683.492 & -1.906 & -6.050 & -1.194 & \textbf{6.072} \\ \bottomrule
        &               &           &            &                &        \\   \\ \toprule
        Results & \textbf{6} & 7 & 8 & 9 & 10 \\ \midrule
        $1$ &   & 17.746 & -0.827 & -14.409 &  \\
        $t$ &   & -7.941 & 1.638 & 6.120 & \\
        $x$  & \textbf{0.195} & -7.443 & 4.711 & 4.546 & -0.163 \\
        $y$ & \textbf{1.000} & -14.308 & -0.968 & 13.091 & -0.837 \\
        $\dot x$ & \textbf{-0.194} & 3.532 & -4.049 & -1.707 & 0.163 \\
        $\dot y$ & \textbf{1} & -12.396 & -1.279 & 11.707 & -0.837  \\
        $tt$ &  & 1 & -0.384 & -0.726 \\ 
        $tx$ &  & -0.206 & 1 & -0.207  \\ 
        $ty$ &   & -0.926 & -0.488 & 1  \\
        $xy$ & \textbf{-1.194} & 16.137 & -3.063 & -13.353 & 1 \\  \bottomrule
    \end{tabular}
\end{table}

\clearpage

\section{Partial differential equations}\label{appendix:pde}
We further explore the ability of LES-SINDy to identify governing equations for three well-known PDEs: the Convection-Diffusion equation, the Burgers equation, and the KS equation. These experiments aim to assess LES-SINDy’s performance in identifying PDEs with varying levels of complexity, which range from linear diffusion processes to nonlinear and chaotic systems.

\subsection{Convection-diffusion equation}
The Convection-Diffusion equation describes the transport of a scalar quantity (such as heat or mass) in a medium where both diffusion and advection processes occur. This experiment tests LES-SINDy’s ability to capture linear processes involving both first-order and second-order derivatives. The results obtained from LES-SINDy for the convection-diffusion equation are presented here.

\begin{table}[!htbp]
\centering\caption{Results of LES-SINDy applied to convection-diffusion equation with varying levels of noise.}\label{tab:ad_nosie}
\begin{tabular}{cc|ccccccccccc}
\toprule
\multicolumn{2}{c|}{\multirow{2}{*}{Results}}     & \multicolumn{8}{c}{candidate functions} \\ 
\multicolumn{2}{c|}{}      & {$u_t$}  & {$u$} & {$u_x$} & {$u_{xx}$} & {$u_{xxx}$} & {$uu_x$} & $uu_{xx}$ & {$uu_{xxx}$} \\ \midrule
\multirow{32}{*}{\rotatebox[origin=c]{90}{noise level}} 
              & truth & 1.0   &       & -1.0   & 1.0    \\ \midrule
 & clean & 1.000 &  & -1.000 & 1.000 &  &  &  &  \\
 & 1\%   & 1.000 &  & -1.000 & 1.000 &  &  &  &  \\
 & 2\%   & 1.000 &  & -1.000 & 1.000 &  &  &  &  \\
 & 3\%   & 1.000 &  & -0.998 & 1.000 &  &  &  &  \\
 & 4\%   & 1.000 &  & -0.996 & 0.999 &  &  &  &  \\
 & 5\%   & 1.000 &  & -0.995 & 0.999 &  &  &  &  \\
 & 6\%   & 1.000 &  & -0.993 & 0.998 &  &  &  &  \\
 & 7\%   & 1.000 &  & -0.991 & 0.998 &  &  &  &  \\
 & 8\%   & 1.000 &  & -0.990 & 0.998 &  &  &  &  \\
 & 9\%   & 1.000 &  & -0.988 & 0.997 &  &  &  &  \\
 & 10\%  & 1.000 &  & -0.986 & 0.996 &  &  &  &  \\
 & 11\%  & 1.000 &  & -0.984 & 0.996 &  &  &  &  \\
 & 12\%  & 1.000 &  & -0.983 & 0.995 &  &  &  &  \\
 & 13\%  & 1.000 &  & -0.981 & 0.995 &  &  &  &  \\
 & 14\%  & 1.000 &  & -0.979 & 0.994 &  &  &  &  \\
 & 15\%  & 1.000 &  & -0.978 & 0.994 &  &  &  &  \\
 & 16\%  & 1.000 &  & -0.976 & 0.993 &  &  &  &  \\
 & 17\%  & 1.000 &  & -0.974 & 0.992 &  &  &  &  \\
 & 18\%  & 1.000 &  & -0.972 & 0.992 &  &  &  &  \\
 & 19\%  & 1.000 &  & -0.971 & 0.991 &  &  &  &  \\
 & 20\%  & 1.000 &  & -0.969 & 0.990 &  &  &  &  \\
 & 21\%  & 1.000 &  & -0.967 & 0.990 &  &  &  &  \\
 & 22\%  & 1.000 &  & -0.965 & 0.989 &  &  &  &  \\
 & 23\%  & 1.000 &  & -0.964 & 0.988 &  &  &  &  \\
 & 24\%  & 1.000 &  & -0.962 & 0.988 &  &  &  &  \\
 & 25\%  & 1.000 &  & -0.960 & 0.987 &  &  &  &  \\
 & 26\%  & 1.000 &  & -0.958 & 0.986 &  &  &  &  \\
 & 27\%  & 1.000 &  & -0.957 & 0.985 &  &  &  &  \\
 & 28\%  & 1.000 &  & -0.955 & 0.985 &  &  &  &  \\
 & 29\%  & 1.000 &  & -0.953 & 0.984 &  &  &  &  \\
 & 30\%  & 1.000 &  & -0.951 & 0.983 &  &  &  &  \\  \bottomrule
\end{tabular}
\end{table}

\clearpage

\subsection{Burgers equation}
The Burgers equation is a fundamental nonlinear PDE that models various physical processes such as fluid dynamics and traffic flow. This equation is notable for its ability to develop shock waves, which makes it a challenging test for identification methods. We present the results from LES-SINDy, demonstrating its ability to capture the nonlinear and diffusive dynamics of the Burgers equation.

\begin{table}[!htbp]
\centering\caption{Results of LES-SINDy applied to Burgers equation with varying levels of noise.}\label{tab:Burgers_noise}
\begin{tabular}{cc|ccccccccccc}
\toprule
\multicolumn{2}{c|}{\multirow{2}{*}{Results}}      & \multicolumn{8}{c}{candidate functions} \\ 
\multicolumn{2}{c|}{}     & {$u_t$}  & {$u$} & $u_x$ & {$u_{xx}$} & {$u_{xxx}$} & {$uu_x$} & $uu_{xx}$ & {$uu_{xxx}$} \\ \midrule
\multirow{22}{*}{\rotatebox[origin=c]{90}{noise level}} 
              & truth & 1.0   &       &        & -0.50 &        & 1.0    \\ \midrule
 & clean & 1.000 &       &  & -0.501 &  & 1.057 &       \\
 & 1\%   & 1.000 &       &  & -0.495 &  & 1.000 &       \\
 & 2\%   & 1.000 &       &  & -0.494 &  & 1.007 & \\
 & 3\%   & 1.000 &       &  & -0.493 &  & 1.000 & \\
 & 4\%   & 1.000 &       &  & -0.492 &  & 0.997 & \\
 & 5\%   & 1.000 &       &  & -0.492 &  & 0.996 & \\
 & 6\%   & 1.000 &       &  & -0.502 &  & 1.057 &       \\
 & 7\%   & 1.000 &       &  & -0.502 &  & 1.057 &       \\
 & 8\%   & 1.000 &       &  & -0.502 &  & 1.057 &       \\
 & 9\%   & 1.000 &       &  & -0.502 &  & 1.058 &       \\
 & 10\%  & 1.000 &       &  & -0.503 &  & 1.058 &       \\
 & 11\%  & 1.000 &       &  & -0.503 &  & 1.059 &       \\
 & 12\%  & 1.000 &       &  & -0.503 &  & 1.060 &       \\
 & 13\%  & 1.000 &       &  & -0.503 &  & 1.061 &       \\
 & 14\%  & 1.000 &       &  & -0.503 &  & 1.062 &       \\
 & 15\%  & 1.000 &       &  & -0.503 &  & 1.063 &       \\
 & 16\%  & 1.000 &       &  & -0.503 &  & 1.065 &       \\
 & 17\%  & 1.000 &       &  & -0.504 &  & 1.066 &       \\
 & 18\%  & 1.000 &       &  & -0.504 &  & 1.067 &       \\
 & 19\%  & 1.000 &       &  & -0.504 &  & 1.069 &       \\
 & 20\%  & 1.000 &       &  & -0.504 &  & 1.070 &       \\ \bottomrule
\end{tabular}
\end{table}

\clearpage
\subsection{Kuramoto–Sivashinsky equation}
The KS equation is a higher-order nonlinear PDE that models instabilities in fluid flows and other dissipative systems. The KS equation is known for its chaotic behavior and complex spatiotemporal dynamics, which makes it a difficult test for model identification algorithms. The results from LES-SINDy for this equation are provided below, which highlights its ability to identify such complex systems.

\begin{table}[!htbp]
\centering\caption{Results of LES-SINDy applied to the KS equation with varying levels of noise.}\label{tab:ks_noise}
\begin{tabular}{cc|ccccccccccc}
\toprule
\multicolumn{2}{c|}{\multirow{2}{*}{Results}}           & \multicolumn{10}{c}{candidate functions} \\ 
\multicolumn{2}{c|}{}     & {$u_t$}  & {$u$} & $u_x$ & {$u_{xx}$} & $u_{xxx}$  & {$u_{xxxx}$} & {$uu_x$} & $uu_{xx}$ & $uu_{xxx}$ & $uu_{xxxx}$ \\ \midrule
\multirow{22}{*}{\rotatebox[origin=c]{90}{noise level}} 
              & truth & 1.0   &       &        &  1.0  &        & 1.0  & 1.0 \\ \midrule
       & clean & 1.000 &  &  & 0.986 &  & 0.995 & 0.987 &  &  &  \\
       & 1\%   & 1.000 &  &  & 0.988 &  & 0.997 & 0.988 &  &  &  \\
       & 2\%   & 1.000 &  &  & 0.989 &  & 0.998 & 0.989 &  &  &  \\
       & 3\%   & 1.000 &  &  & 0.990 &  & 0.999 & 0.989 &  &  &  \\
       & 4\%   & 1.000 &  &  & 0.990 &  & 1.000 & 0.989 &  &  &  \\
       & 5\%   & 1.000 &  &  & 0.990 &  & 1.000 & 0.989 &  &  &  \\
       & 6\%   & 1.000 &  &  & 0.990 &  & 1.000 & 0.988 &  &  &  \\
       & 7\%   & 1.000 &  &  & 0.989 &  & 1.000 & 0.987 &  &  &  \\
       & 8\%   & 1.000 &  &  & 0.988 &  & 0.999 & 0.986 &  &  &  \\
       & 9\%   & 1.000 &  &  & 0.986 &  & 0.998 & 0.984 &  &  &  \\
       & 10\%  & 1.000 &  &  & 0.983 &  & 0.996 & 0.982 &  &  &  \\
       & 11\%   & 1.000 &  &  & 0.981 &  & 0.994 & 0.979 &  &  &  \\
       & 12\%   & 1.000 &  &  & 0.978 &  & 0.991 & 0.976 &  &  &  \\
       & 13\%   & 1.000 &  &  & 0.974 &  & 0.988 & 0.973 &  &  &  \\
       & 14\%   & 1.000 &  &  & 0.970 &  & 0.985 & 0.969 &  &  &  \\
       & 15\%   & 1.000 &  &  & 0.965 &  & 0.981 & 0.965 &  &  &  \\
       & 16\%   & 1.000 &  &  & 0.961 &  & 0.977 & 0.960 &  &  &  \\
       & 17\%   & 1.000 &  &  & 0.955 &  & 0.973 & 0.956 &  &  &  \\
       & 18\%   & 1.000 &  &  & 0.949 &  & 0.968 & 0.950 &  &  &  \\
       & 19\%   & 1.000 &  &  & 0.943 &  & 0.962 & 0.945 &  &  &  \\
       & 20\%   & 1.000 &  &  & 0.937 &  & 0.957 & 0.939 &  &  &  \\ \bottomrule
\end{tabular}
\end{table}

\end{appendices}

\end{document}